\documentstyle[%showkeys,
11pt,amsfonts]{article}

\newcommand {\rel} {{\mathbb R}}
\newcommand {\com} {{\mathbb C}}
\newcommand {\nat} {{\mathbb N}}

\newcommand {\ganz} {{\mathbb Z}}

\newtheorem{proposition}{Proposition}[section]
\newtheorem{theorem}{Theorem}[section]

\renewcommand{\theequation}{\mbox{\arabic{section}.\arabic{equation}}}

\newcommand{\nc}{\newcommand}

\newcommand{\bcite}[1] {\cite{#1}}

\newcommand {\pr} {\bf}

\newcommand{\proof} {   \begin{flushright}
                        ///
                        \end{flushright}
                }

\newcommand{\defin} { \hspace*{\fill} $\Box$ }

\newcommand{\foot}[1] { }
\renewcommand{\foot}[1] { \footnote{#1} }

\def    \curro  {{\bf [ \hskip -1.5pt [ }}
\def    \currc  {{\bf ] \hskip -1.5pt ] }}
\newcommand     {\curr}[1]      {{ \curro #1 \currc }}

\def    \mean   {{ \vec {\bf H} }}

\def    \Ht	{{ {\cal H}^2 }}

\def    \Lt	{{ {\cal L}^2 }}

\def	\d	{{ \ {\rm d} }}

\def    \thez   {{ \theta_0 }}
\def    \theo   {{ \theta }}

\def    \W      {{ {\mathcal W} }}
\nc{\energ}[1]	{{ e({#1}) }}

\def    \Cn     {{ C_{n} }}

\def	\spt	{{ spt\ }}

\def	\poin	{{ poin }}

\def	\geu	{{ g_{euc} }}

\nc{\gpo}[1]	{{ g_{\poin}^{#1} }}
\nc{\gpu}[1]	{{ g_{\poin,{#1}} }}

\nc{\gpmo}[1]	{{ g_{\poin,m}^{#1} }}
\nc{\gpmu}[1]	{{ g_{\poin,m,{#1}} }}
\nc{\tgpmo}[1]	{{ \tilde g_{\poin,m}^{#1} }}

\nc{\gplo}[1]	{{ g_{\poin,\lambda}^{#1} }}

\nc{\conf}[1]	{{ [{#1}] }}
\def	\nor	{{ \cal N }}

\def    \diff	     {{ \stackrel{\approx}{\longrightarrow} }}

\nc{\doub}[1]{{ \ddot{#1} }}
\nc{\dd}{ \begin{displaymath} }
\nc{\df}{ \end{displaymath} }
\nc{\dcd}{ \begin{displaymath} \begin{array}{c}}
\nc{\dcf}{ \end{array} \end{displaymath} }
\nc{\ee}{ \begin{equation} }
\nc{\ef}{ \end{equation} }
\nc{\ad}{ \begin{array}{c} }
\nc{\af}{ \end{array} }

\begin{document}

%\begin{flushright}
%not for distribution
%\end{flushright}
\begin{center}
{\huge \bf Optimal rigidity estimates for nearly umbilical surfaces
in arbitrary codimension}
\\ \ \\
Tobias Lamm \\
Fakult\"at f\"ur Mathematik\\
Karlsruher Institut f\"ur Technologie (KIT)\\
Kaiserstra\ss e 89-93,
D-76133 Karlsruhe, Germany, \\
email: tobias.lamm@kit.edu \\
\ \\
Reiner Michael Sch\"atzle \\
Fachbereich Mathematik der
Eberhard-Karls-Universit\"at T\"ubingen, \\
Auf der Morgenstelle 10,
%Geb\"aude C, 5 A 40,
D-72076 T\"ubingen, Germany, \\
email: schaetz@mail.mathematik.uni-tuebingen.de \\
 
\end{center}
\vspace{1cm}

\begin{quote}

{\bf Abstract:}
In \bcite{del.muell},
DeLellis and M\"uller proved a quantitative version of Codazzi's theorem,
namely for a smooth embedded surface $\ \Sigma \subseteq \rel^3\ $
with area normalized to $\ \Ht(\Sigma) = 4 \pi\ $,
it was shown that
$\ \parallel A_\Sigma - id \parallel_{L^2(\Sigma)}
\leq C \parallel A^0_\Sigma \parallel_{L^2(\Sigma)}\ $,
and building on this,
closeness of $\ \Sigma\ $ to a round sphere in $\ W^{2,2}\ $ was established,
when $\ \parallel A^0_\Sigma \parallel_{L^2(\Sigma)}\ $ is small.
This was supplemented in \bcite{del.muell.2}
by giving a conformal parametrization $\ S^2 \diff \Sigma\ $
with small conformal factor in $\ L^\infty\ $,
again when $\ \parallel A^0_\Sigma \parallel_{L^2(\Sigma)}\ $ is small.
In this article, we extend these results to arbitrary codimension.
In contrast to \bcite{del.muell},
our argument is not based on the equation of Mainardi-Codazzi,
but instead uses the monotonicity formula for varifolds.
\end{quote}

\vspace{1cm}

%\tableofcontents

%%%%%

\setcounter{equation}{0}

\section{Introduction} \label{intro}

A classical theorem of Codazzi in differential geometry states
that all umbilical connected surfaces,
that is their tracefree second fundamental form $\ A^0 \equiv 0\ $ vanishes,
are pieces of a round sphere or a plane.
A quantitive version of this theorem was given
by DeLellis and M\"uller in \bcite{del.muell} in codimension 1,
namely for a smooth embedded surface $\ \Sigma \subseteq \rel^3\ $
with area normalized to $\ \Ht(\Sigma) = 4 \pi\ $,
it was shown that the second fundamental form $A_\Sigma $ of $\Sigma$ satisfies
\begin{equation} \label{intro.funda}
	\parallel A_\Sigma - id \parallel_{L^2(\Sigma)}
	\leq C \parallel A^0_\Sigma \parallel_{L^2(\Sigma)},
\end{equation}
and if
\begin{equation} \label{intro.funda-0}
	\parallel A^0_\Sigma \parallel_{L^2(\Sigma)}^2 < 4 \pi,
\end{equation}
there exists a conformal parametrization
$\ \psi: S^2 \diff \Sigma \mbox{ and } c_\Sigma \in \rel^3\ $ with
\begin{equation} \label{intro.round}
	\parallel \psi - (c_\Sigma + id_{S^2}) \parallel_{W^{2,2}(S^2)}
	\leq C \parallel A^0_\Sigma \parallel_{L^2(\Sigma)}.
\end{equation}
Here and in the following $\ g, \mean, A^0 = A - \frac{1}{2} \mean g\ $
denote the pull-back metric, the mean curvature
and the tracefree second fundamental form of $\ \Sigma \ $.
 
In \bcite{del.muell.2}, the estimate (\ref{intro.round}) was supplemented
with a $\ L^\infty-$bound
on the conformal factor of the conformal
pull-back metric $\ g = \psi^* \geu = h^2 g_{S^2}\ $
\begin{equation} \label{intro.conf}
	\parallel h - 1 \parallel_{L^\infty(S^2)}
	\leq C \parallel A^0_\Sigma \parallel_{L^2(\Sigma)}.
\end{equation}
Nguyen and the first author gave in \bcite{lamm.nguy.conf}
an extension to arbitrary codimension
in the form that if $\ \parallel A^0_\Sigma \parallel_{L^2(\Sigma)}^2
\rightarrow 0 \mbox{ and } \Ht(\Sigma) = 4 \pi\ $,
then for some conformal parametrization $\ \psi: S^2 \diff \Sigma\ $
after an appropriate translation and rotation
\begin{displaymath}
	\parallel \psi - id_{S^2} \parallel_{W^{2,2}(S^2)} \rightarrow 0.
\end{displaymath}
In this article, we extend (\ref{intro.funda}),
(\ref{intro.round}) and (\ref{intro.conf})
after an appropriate translation and rotation to any codimension
and improve the bound in (\ref{intro.funda-0}).

\begin{theorem} \label{intro.theorem-funda}

Let $\ \Sigma \subseteq \rel^n\ $ be a smoothly embedded
closed connected surface with $\ \Ht(\Sigma) = 4 \pi\ $.
Then there exists a measurable unit normal vector field
$\ \nor \mbox{ on } \Sigma\ $ with
\begin{displaymath}
	\parallel A_\Sigma - \nor g \parallel_{L^2(\Sigma)}
	\leq \Cn \parallel A^0_\Sigma \parallel_{L^2(\Sigma)}
\end{displaymath}
and
\begin{displaymath}
	\parallel K_\Sigma - 1 \parallel_{L^1(\Sigma)}
	\leq \Cn \parallel A^0_\Sigma \parallel_{L^2(\Sigma)}^2,
\end{displaymath}
where $K_\Sigma$ denotes the Gauss curvature of $\Sigma$.
\defin
\end{theorem}
{\large \bf Remark:} \\
Using the approximation technique of Schoen and Uhlenbeck
in \bcite{schoe.uhl.harmb} \S 4 Proposition,
this extends to uniformly conformal $W^{2,2}$-immersions.
\defin

\begin{theorem} \label{intro.theorem-round}

Let $\ \Sigma \subseteq \rel^n\ $ be a smoothly embedded surface
of sphere type $\ \Sigma \cong S^2\ $
with $\ \Ht(\Sigma) = 4 \pi\ $ and
\begin{displaymath}
	\parallel A^0_\Sigma \parallel_{L^2(\Sigma)}^2
	< 2 \energ{n}
	= \left\{
	\begin{array}{cl}
		8 \pi & \mbox{for } n = 3, \\
		16 \pi / 3 & \mbox{for } n = 4, \\
		4 \pi &  \mbox{for } n \geq 5, \\
	\end{array}
	\right.
\end{displaymath}
where $\ \energ{n} = e_n\ $ was defined in \bcite{schae.comp-will12} (1.2).

Then there exists a conformal parametrization $\ f: S^2 \diff \Sigma\ $
with pull-back metric $\ g = f^* \geu = e^{2u}  g_{S^2}\ $,
such that after an appropriate translation and rotation
and with $\ S^2 := \partial B_1(0)
\cap span \{ e_1 , e_2 , e_3 \} \subseteq \rel^n\ $
\begin{displaymath}
	\parallel f - id_{S^2} \parallel_{W^{2,2}(S^2)}
	+ \parallel u \parallel_{L^\infty(S^2)}
	\leq C(n,\tau) \parallel A^0_\Sigma \parallel_{L^2(\Sigma)},
\end{displaymath}
where $\ \tau := 2 \energ{n}
- \parallel A^0_\Sigma \parallel_{L^2(\Sigma)}^2 > 0\ $.
\defin
\end{theorem}
{\large \bf Remark:} \\
The bounds on the right-hand side
of the assumption cannot be improved beyond $\ 8 \pi\ $,
as two spheres connected by a small part of a catenoid show.
In particular, the assumption is optimal for $\ n = 3\ $.
\defin
\\ \ \\
In \bcite{del.muell}, an important ingredient of the proof
were the analytical Hardy-space estimates in \bcite{muell.sver}
and the equation of Mainardi-Codazzi, which was used to establish
\begin{equation} \label{intro.aux}
	\Big\|\ \Big| \frac{1}{2} \mean_\Sigma \Big| - 1 \Big\|_{L^2(\Sigma)}
	\leq C \parallel A^0_\Sigma \parallel_{L^2(\Sigma)}
\end{equation}
to obtain (\ref{intro.funda}).

Instead of using the equation of Mainardi-Codazzi,
our proof relies on the observation that the Willmore functional,
which is the square integral of the mean curvature
$\ \mean_\Sigma \mbox{ of } \Sigma\ $
multiplied by a factor $\ 1/4\ $ that is
\begin{displaymath}
	\W(\Sigma) := \frac{1}{4} \int \limits_\Sigma |\mean_\Sigma|^2 \d \Ht,
\end{displaymath}
satisfies for closed $\ \Sigma\ $
\begin{displaymath}
	\W(\Sigma) \geq 4 \pi,
\end{displaymath}
and the global minimizers are the round spheres,
see \bcite{will.conj} in $\ \rel^3\ $.
In general codimension,
the inequality is a consequence of the Li-Yau inequality,
see \bcite{li.yau}.

By the Gau\ss\ equations and the Gau\ss-Bonnet theorem, we have
\begin{equation} \label{intro.gauss}
	\W(\Sigma)
	= \frac{1}{4} \int \limits_{\Sigma}
	|A_\Sigma|^2 \d \Ht + \pi \chi(\Sigma)
	= \frac{1}{2} \int \limits_{\Sigma}
	|A^0_\Sigma|^2 \d \Ht + 2 \pi \chi(\Sigma),
\end{equation}
and see for $\ \Sigma \cong S^2\ $ of sphere type that
\begin{displaymath}
	\W(\Sigma) = 4 \pi
	+ \frac{1}{2} \parallel A^0_\Sigma \parallel_{L^2(\Sigma)}^2,
\end{displaymath}
which shows the equivalence of the smallness assumption of
$\ \parallel A^0_\Sigma \parallel_{L^2(\Sigma)}\ $
and the Willmore energy $\ \W(\Sigma)\ $
being close to the absolute minimum $\ 4 \pi\ $,
when $\ \Sigma \cong S^2\ $.
In any case, smallness of $\ \parallel A^0_\Sigma \parallel_{L^2(\Sigma)}\ $
implies $\ \W(\Sigma) \approx 4 \pi\ $,
as $\ \chi(\Sigma) \leq 2\ $,
and $\ \Sigma \cong S^2\ $ is of sphere type.

In \S \ref{glob}, we give an argument using the monotonicity formula,
see \bcite{sim.will},
that the global minimizers of the Willmore functional are the round spheres,
which even works in the non-smooth case.
More precisely, the non-negative term in the monotonicity formula yields
for $\ \Ht(\Sigma) = 4 \pi\ $ and
\begin{displaymath}
	\W(\Sigma) \leq 4 \pi + \delta^2
\end{displaymath}
that
\begin{displaymath}
	\int \limits_\Sigma \Big| \frac{1}{4} \mean_\Sigma(y)
	+ \frac{(y - x)^{\perp_y}}{|y-x|^2} \Big|^2
	\d \Ht(y) \leq \delta^2/4
	\quad \mbox{for all } x \in \Sigma.
\end{displaymath}
where $\ ^{\perp_y}\ $ denotes the orthogonal projection
onto $\ N_y \Sigma\ $.
Then by Fubini's theorem
\begin{displaymath}
	\int \limits_\Sigma \Big| \frac{1}{4} \mean_\Sigma(y)
	+ \frac{(y - x)^{\perp_y}}{|y-x|^2} \Big|^2
	\d \Ht(x) \leq \delta^2/4
	\quad \mbox{for some } y \in \Sigma.
\end{displaymath}
Observing for the two-dimensional round sphere
$\ S\ $ tangent to $\ \Sigma \mbox{ at } y\ $
and with radius $\ 2 / |\mean_\Sigma(y)|
\mbox{ in } T_y \Sigma + span \{ \mean_\Sigma(y) \}\ $,
if $\ \mean_\Sigma(y) \neq 0\ $, that
\begin{displaymath}
	2 d(x,S) \leq \Big| \mean_\Sigma(y)
	+ 4 \frac{(y-x)^\perp}{|y-x|^2} \Big|\ |y-x|^2,
\end{displaymath}
we obtain
\begin{displaymath}
	\int \limits_\Sigma d(x,S)^2 \d \Ht(x)
	\leq C (diam\ \Sigma)^4 \delta^2,
\end{displaymath}
and similarly for $\ S\ $ replaced by the tangent plane $\ T_y \Sigma\ $,
if $\ \mean_\Sigma(y) = 0\ $.
Examining all possible cases,
we conclude that $\ \Sigma\ $ is close to a two-dimensional round sphere of radius 1
in the sense that after an appropriate translation and rotation
\begin{displaymath}
	d_H(\Sigma , S^2) \leq \Cn \sqrt{\delta},
\end{displaymath}
where $\ d_H\ $ denotes the Hausdorff distance,
\begin{equation} \label{intro.mean}
	\parallel \mean_\Sigma + 2 id_\Sigma \parallel_{L^2(\Sigma)}
	\leq \Cn \delta,
\end{equation}
\begin{displaymath}
	\parallel |id_\Sigma| - 1 \parallel_{L^2(\Sigma)}
	\leq \Cn \delta,
\end{displaymath}
see Propositions \ref{glob.unit-sphere} and \ref{glob-smooth.sphere}.
This already yields (\ref{intro.aux})
for $\ \delta\ $ small enough depending on $\ n\ $.
Then the first estimate in Theorem \ref{intro.theorem-funda}
immediately follows,
and the second estimate follows
from a general estimate on the Gau\ss-curvature,
see Proposition \ref{glob-smooth.gauss}.
For $\ \parallel A^0_\Sigma \parallel_{L^2(\Sigma)}
\geq \delta > 0\ $,
we see from (\ref{intro.gauss}) that
\begin{equation} \label{intro.gauss-impl}
	\parallel A_\Sigma \parallel_{L^2(\Sigma)}
	\leq C(\delta) \parallel A^0_\Sigma \parallel_{L^2(\Sigma)},
\end{equation}
since $\ \chi(\Sigma) \leq 2\ $,
and Theorem \ref{intro.theorem-funda} is immediate,
see the end of \S \ref{glob}.

In \S \ref{conf},
we get the conformal parametrization in Theorem \ref{intro.theorem-round},
first by parametrizing $\ \Sigma\ $ by the uniformization theorem.
Then inverting $\ \Sigma\ $ and parametrizing over $\ \com\ $,
we can estimate the conformal factor on the plane
with the Hardy space estimates in \bcite{muell.sver}.
In particular, we use the bi-Lipschitz estimates in \bcite{muell.sver}
Theorem 4.3.1 together with the improvement of energy constants
in \bcite{kuw.schae.will6} Theorem 6.1 for $\ n = 3,4\ $,
and \bcite{schae.comp-will12} Theorem 5.1. for every $n\ge 3$.
Applying a dilation and a translation in the plane,
the conformal factor on the sphere is bounded as well.

Our arguments in \S \ref{conf} differ quite a bit from the corresponding ones in \bcite{del.muell} (see Proposition 3.2 therein) since we have not been able to directly extend the uniform bounds for the conformal factor to arbitrary codimensions. Moreover, since we use the bi-Lipschitz estimates from \bcite{muell.sver} we have been able to simplify several of the arguments in \bcite{del.muell}, \bcite{del.muell.2}, in particular the smallness of the conformal factor in $L^\infty$ is more or less a direct consequence of our new approach and these estimates (see Proposition \ref{conf.small}). 

In \S \ref{round}, we observe from (\ref{intro.mean}) and $\ \Delta_g f = \mean_f\ $
that $\ f\ $ nearly lies in the kernel of $\ \Delta_{S^2} + 2\ $,
which consists of precisely the linear functions.
Then following the estimate in \bcite{del.muell} \S 6,
we get the $\ W^{2,2}-$bound in Theorem \ref{intro.theorem-round}.
The smallness of the conformal factor in Theorem \ref{intro.theorem-round}
is obtained by the inversion from \S \ref{conf}
when knowing that $\ \Sigma \mbox{ is } L^\infty-$close to a round sphere,
which is implied by the $\ W^{2,2}-$bound and the Sobolev embedding. In this section parts of our arguments are direct modifications of the corresponding results in \bcite{del.muell}, see e.g. Proposition 4.1, whereas we also introduced a streamlined argument for the actual $W^{2,2}$-closedness of the surfaces in Proposition 4.2. In particular, we circumvented the use of the Cartan formalism.

We note that the optimal rigidity estimates of DeLellis and M\"uller were crucial ingredients in the construction of foliations of asymptotically flat resp. asymptotically hyperbolic $3$-manifolds by surfaces of prescribed mean curvature \bcite{metzger07}, by surfaces of Willmore type \bcite{lam.met.sch} resp. by surfaces of constant mean curvature \bcite{nev.tia}. Moreover, they were used in order to study spherical critical points of $\W$ with prescribed area in Riemannian $3$-manifolds, see \bcite{lam.met}, \bcite{lam.met2} and \bcite{muell.roeger}. Finally, they were used in a result of R\"oger and the second author \bcite{roeger.schae.will} were an estimate for the isoperimetric deficit in terms of the Willmore deficit was derived. 
We anticipate that Theorem \ref{intro.theorem-funda} and Theorem \ref{intro.theorem-round} will turn out to be crucial ingredients for extensions of the above results to higher codimensions.

%%%%%
%%%%%

\setcounter{equation}{0}

\section{Global Willmore minimizers} \label{glob}

It is known
that the global minimizers of the Willmore energy
are the round spheres,
see \bcite{will.conj} in $\ \rel^3\ $.
Here we give an argument which works in any codimension
without assuming regularity
by using the monotonicity formula
developed by Simon in \bcite{sim.will}
and continued in \bcite{kuw.schae.will3}.
For the notions in geometric measure theory,
we refer to \bcite{sim}.

\begin{proposition} \label{glob.mini}

Let $\ \mu \neq 0\ $ be an integral $\ 2-$varifold
with square integrable weak mean curvature
and compact support.
Then
\begin{displaymath}
	\W(\mu) \geq 4 \pi
\end{displaymath}
and in case of equality,
$\ \mu\ $ is a single round sphere.
\end{proposition}
{\pr Proof:} \\
The above inequality was already obtained
in \bcite{kuw.schae.will3} (A.18).
In case of equality,
we consider any  $\ x \in spt\ \mu \neq \emptyset\ $
and see by monotonicity in \bcite{kuw.schae.will3} (A.3) of
\begin{equation} \label{glob.mini.gamma}
	\gamma(\varrho) :=  \varrho^{-2} \mu(B_\varrho(x))
	+ \frac{1}{16} \int \limits_{B_\varrho(x)}
	|\mean_\mu|^2 \d \mu
	+ \frac{1}{2} \int \limits_{B_\varrho(x)}
	\varrho^{-2} (y - x)\ \mean_\mu(y) \d \mu(y)
\end{equation}
defined in \bcite{kuw.schae.will3} (A.4),
and by \bcite{kuw.schae.will3} (A.7), (A.10), (A.14) that
\begin{displaymath}
	\pi \leq \omega_2 \theta^2(\mu,x)
	= \lim \limits_{\varrho \rightarrow 0} \gamma(\varrho)
	\leq \gamma(r)
	\leq \lim \limits_{\varrho \rightarrow \infty} \gamma(\varrho)
	= \frac{1}{4} \W(\mu) = \pi 
	\quad \forall r,
\end{displaymath}
hence $\ \gamma \equiv \pi\ $ is constant.
Then
\begin{equation} \label{glob.mini.dense}
	\theta^2(\mu) = 1
	\quad \mbox{on } spt\ \mu
\end{equation}
and by \bcite{kuw.schae.will3} (A.3)
\begin{displaymath}
	\mean_\mu(y) + 4 \frac{(y - x)^{\perp_y}}{|y - x|^2} = 0
	\quad \mbox{for } \mu-\mbox{almost all } y \in spt\ \mu,
\end{displaymath}
where $\ ^{\perp_y}\ $ denotes the orthogonal projection
onto $\ N_y \mu\ $.
In particular
\begin{equation} \label{glob.mini.perp}
	\mean_\mu(y) \perp T_y \mu
	\quad \mbox{for } \mu-\mbox{almost all } y \in spt\ \mu.
\end{equation}
By Fubini's theorem, we get for $\ \mu-\mbox{almost all } y\ $ that
\begin{displaymath}
	\mean_\mu(y) + 4 \frac{(y - x)^{\perp_y}}{|y - x|^2} = 0
	\quad \mbox{for } \mu-\mbox{almost all } x \in spt\ \mu.
\end{displaymath}
We choose any such $\ y \in spt\ \mu\ $,
in particular $\ T_y \mu\ $ exists.
If $\ \mean_\mu(y) = 0\ $,
then $\ spt\ \mu \subseteq y + T_y \mu\ $,
in particular $\ T_x \mu = T_y \mu\ $
and by (\ref{glob.mini.perp}) that
$\ \mean_\mu(x) \perp T_y \mu
\mbox{ for } \mu-\mbox{almost all } x \in spt\ \mu\ $.
Then $\ \mu\ $ is stationary in $\ y + T_y \mu\ $
in the sense of \bcite{sim} 16.4 or 41.2 (3),
and by constancy theorem, see \bcite{sim} Theorem 41.1,
we get $\ \mu = \theta \Ht \lfloor (y + T_y \mu)
\mbox{ for some constant } \theta > 0\ $.
This contradicts the compactness of $\ \spt\ \mu\ $.
Hence $\ \mean_\mu(y) \neq 0\ $,
and we may further assume that
$\ \mean_\mu(y) \perp T_y \mu\ $ by (\ref{glob.mini.perp}).
%see \bcite{brakke} Theorem 5.8.
To abbreviate notation,
we assume after rotation, scaling and translation
that $\ y = 0 \mbox{ and } T_0 \mu = span \{ e_1 , e_2 \},
N_0 \mu = span \{ e_3, \ldots, e_n \},
\mean_\mu(0) = 2 e_3\ $
and write $\ ^\perp \mbox{ for } \ ^{\perp_y}\ $.
We firstly get from above for $\ j = 4, \ldots, n\ $, that
\begin{displaymath}
	0 = \langle \mean_\mu(0) , e_j \rangle
	= - 4 \Big\langle \frac{-x^\perp}{|x|^2} , e_j \Big\rangle =
\end{displaymath}
\begin{displaymath}
	= 4 \Big\langle \frac{x^\perp}{|x|^2} , e_j \Big\rangle
	= 4 x_j / |x|^2
	\quad \mbox{for } \mu-\mbox{almost all } x \neq 0 \in spt\ \mu,
%	j = 4, \ldots, n,
\end{displaymath}
hence $\ spt\ \mu \subseteq \rel^3\ $.
Next for $\ j = 3\ $
\begin{displaymath}
	2 = \langle \mean_\mu(0) , e_3 \rangle
	= - 4 \Big\langle \frac{-x^\perp}{|x|^2} , e_3 \Big\rangle
%	= 4 \Big\langle \frac{x^\perp}{|x|^2} , e_3 \Big\rangle
	= 4 x_3 / |x|^2
	\quad \mbox{for } \mu-\mbox{almost all } x \neq 0 \in spt\ \mu,
\end{displaymath}
hence $\ 2 x_3 = |x|^2\ $ or likewise
\begin{displaymath}
	1 = x_1^2 + x_2^2 + x_3^2 - 2 x_3 + 1
	= x_1^2 + x_2^2 + (1 - x_3)^2
	= |x - e_3|^2
\end{displaymath}
and $\ spt\ \mu \subseteq \partial_1 B(e_3)\ $.
Together we see that $\ spt\ \mu \subseteq
\partial B_1(e_3) \cap \rel^3 =: S \cong S^2\ $,
in particular $\ T_y \mu = T_y S\ $
and by (\ref{glob.mini.perp}) that
$\ \mean_\mu(y) \perp T_y S
\mbox{ for } \mu-\mbox{almost all } y \in spt\ \mu\ $.
Then $\ \mu\ $ is stationary in $\ S\ $
in the sense of \bcite{sim} 16.4 or 41.2 (3),
and again by constancy theorem, see \bcite{sim} Theorem 41.1,
we get $\ \mu = \theta \Ht \lfloor S
\mbox{ for some constant } \theta > 0\ $.
By (\ref{glob.mini.dense}), we see that $\ \theta = 1\ $,
and the proposition is proved.
\proof
Actually we can give a more quantitative version
when the Willmore energy is only close
to the minimal energy of $\ 4 \pi\ $.

\begin{proposition} \label{glob.unit-sphere}

Let $\ \mu \neq 0\ $ be an integral $\ 2-$varifold
with square integrable weak mean curvature,
compact support and
\begin{equation} \label{glob.unit-sphere.measure}
	\mu(\rel^n) = 4 \pi
\end{equation}
and with an appropriate orientation
let $\ T\ $ be an integral current with
underlying measure $\ \mu_T = \mu\ $
and with boundary $\ \partial T = 0\ $.
If
\begin{equation} \label{glob.unit-sphere.energy}
	\W(\mu) \leq 4 \pi + \delta^2,
\end{equation}
then for small $\ \delta\ $ depending on $\ n\ $
there exists a two dimensional round sphere
$\ S_1 \subseteq \rel^n\ $ with radius 1 and
\begin{equation} \label{glob.unit-sphere.sphere}
	\int d(x,S_1)^2 \d \mu(x) \leq \Cn \delta^2,
\end{equation}
\begin{equation} \label{glob.unit-sphere.haus}
	d_H(spt\ \mu , S_1) \leq \Cn \sqrt{\delta},
\end{equation}
where $\ d_H\ $ denotes the Hausdorff distance,
and for $\ S_1\ $ to be centered at the origin
\begin{equation} \label{glob.unit-sphere.mean}
	\int | \mean_\mu(x) + 2 x |^2 \d \mu(x)
	\leq \Cn \delta^2
\end{equation}
for some $\ \Cn < \infty\ $ depending only on the codimension.
\end{proposition}
{\pr Proof:} \\
We may assume that $\ \W(\mu) < 8 \pi \mbox{ for } \delta^2 < 4 \pi\ $
and then see that $\ spt\ \mu\ $ is connected by (\ref{haus.mono.conn}).
By (\ref{glob.unit-sphere.measure}),
we get from \bcite{sim.will} Lemma 1.1
and (\ref{haus.meas-low.esti})
by connectedness of $\ spt\ \mu\ $ for the non-smooth case
and by the density bound in \bcite{kuw.schae.will3} (A.16), i.e.
\begin{equation} \label{glob.unit-sphere.density}
	\varrho^{-2} \mu(B_\varrho)) \leq C
	\quad \mbox{for any } B_\varrho \subseteq \rel^n
\end{equation}
that
\begin{equation} \label{glob.unit-sphere.diam}
	c_0 \leq diam\ spt\ \mu \leq \Cn
\end{equation}
for some $\ c_0 > 0, \Cn < \infty\ $.
Moreover $\ \theta^2(\mu) \leq \W(\mu) / (4 \pi) < 2\ $
by \bcite{kuw.schae.will3} (A.17),
hence as $\ \mu\ $ is integral
\begin{equation} \label{glob.unit-sphere.support}
	\mu = \Ht \lfloor spt\ \mu.
\end{equation}
Again for $\ \gamma\ $ as in (\ref{glob.mini.gamma}),
we see from \bcite{kuw.schae.will3} (A.3) that
\begin{equation} \label{glob.unit-sphere.mono}
	\gamma(\varrho) - \gamma(\sigma)
	= \int \limits_{B_\varrho(x) - B_\sigma(x)}
	\Big| \frac{1}{4} \mean_\mu(y)
	+ \frac{(y - x)^{\perp_y}}{|y-x|^2} \Big|^2
	\d \mu(y) \geq 0
	\quad \mbox{for } 0 < \sigma \leq \varrho < \infty,
\end{equation}
where $\ ^{\perp_y}\ $ denotes the orthogonal projection
onto $\ N_y \mu\ $.
Again by \bcite{kuw.schae.will3} (A.7), (A.10), (A.14),
we get for any $\ x \in spt\ \mu \neq \emptyset\ $
by (\ref{glob.unit-sphere.energy}) that
\begin{displaymath}
	\lim \limits_{\varrho \rightarrow \infty} \gamma(\varrho)
	- \lim \limits_{\sigma \rightarrow 0} \gamma(\sigma)
	= \frac{1}{4} \W(\mu) - \omega_2 \theta(\mu,x)
	\leq \delta^2/4,
\end{displaymath}
hence
\begin{equation} \label{glob.unit-sphere.int}
	\int \Big| \frac{1}{4} \mean_\mu(y)
	+ \frac{(y - x)^{\perp_y}}{|y-x|^2} \Big|^2
	\d \mu(y) \leq \delta^2/4
	\quad \mbox{for all } x \in spt\ \mu.
\end{equation}
Integrating $\ x \mbox{ by } \mu\ $, we get by Fubini's theorem
\begin{displaymath}
	\int \int \Big| \frac{1}{4} \mean_\mu(y)
	+ \frac{(y - x)^{\perp_y}}{|y-x|^2} \Big|^2
	\d \mu(x) \d \mu(y) \leq \mu(\rel^n) \delta^2/4,
\end{displaymath}
hence
\begin{equation} \label{glob.unit-sphere.int-fub}
	\int \Big| \mean_\mu(\xi)
	+ 4 \frac{(\xi - x)^{\perp_\xi}}{|\xi-x|^2} \Big|^2
	\d \mu(x) \leq 4 \delta^2
	\quad \mbox{for some } \xi \in spt\ \mu,
\end{equation}
for which $\ T_\xi \mu \mbox{ and } \mean_\mu(\xi) \in \rel^n\ $ exist
and for which $\ \mean_\mu(\xi) \perp T_\xi \mu\ $
by \bcite{brakke} Theorem 5.8.
To abbreviate notation,
we assume after rotation and translation
that $\ \xi = 0 \mbox{ and } T_0 \mu = span \{ e_1 , e_2 \},
N_0 \mu = span \{ e_3, \ldots, e_n \},
\mean_\mu(0) = 2 \alpha e_3 \mbox{ with } \alpha \geq 0\ $
and write $\ ^\perp \mbox{ for }\ ^{\perp_\xi}\ $.
If $\ \mean_\mu(0) = 0\ $, we have
\begin{displaymath}
	|x|^2\ \Big| \mean_\mu(0) + 4 \frac{(-x)^\perp}{|x|^2} \Big|
	= 4 |x^\perp| = 4 d(x , span \{ e_1 , e_2 \} ),
\end{displaymath}
hence by (\ref{glob.unit-sphere.diam}) and (\ref{glob.unit-sphere.int-fub})
\begin{equation} \label{glob.unit-sphere.dist-plane}
	\int d(x , span \{ e_1 , e_2 \} )^2 \d \mu(x) \leq \Cn \delta^2.
\end{equation}
If $\ \mean_\mu(0), \alpha \neq 0\ $,
and we get as in the previous proof for $\ x \neq 0\ $
\begin{displaymath}
	\Big\langle \mean_\mu(0)
	+ 4 \frac{-x^\perp}{|x|^2} , e_j \Big\rangle
	= -4 \Big\langle \frac{x^\perp}{|x|^2} , e_j \Big\rangle
	= -4 x_j / |x|^2
	\quad \mbox{for } j = 4, \ldots, n,
\end{displaymath}
\begin{displaymath}
	\Big\langle \mean_\mu(0)
	+ 4 \frac{-x^\perp}{|x|^2} , e_3 \Big\rangle
	= 2 \alpha - 4 x_3 / |x|^2
	= 2 \frac{\alpha |x|^2 - 2 x_3}{|x|^2}
	= 2 \frac{|\sqrt{\alpha} x - e_3/\sqrt{\alpha}|^2 - 1/\alpha}{|x|^2}
\end{displaymath}
and
\begin{displaymath}
	|x|^2\ \Big| \Big\langle \mean_\mu(0)
	+ 4 \frac{-x^\perp}{|x|^2} , e_3 \Big\rangle \Big|
	= 2 |\ |\sqrt{\alpha} x - e_3/\sqrt{\alpha}|^2 - 1/\alpha\ | =
\end{displaymath}
\begin{displaymath}
	= 2 |\ |\sqrt{\alpha} x - e_3/\sqrt{\alpha}| + 1/\sqrt{\alpha}\ |
	\cdot |\ |\sqrt{\alpha} x - e_3/\sqrt{\alpha}|
	- 1/\sqrt{\alpha}\ | \geq
\end{displaymath}
\begin{displaymath}
	\geq (2/\sqrt{\alpha}) |\ |\sqrt{\alpha} x - e_3/\sqrt{\alpha}|
	- 1/\sqrt{\alpha}\ |
	= 2 |\ |x - e_3/\alpha| - 1/\alpha\ |.
\end{displaymath}
Putting $\ S_r := \partial B_{1/\alpha}(e_3/\alpha)
\cap span \{ e_1 , e_2 , e_3 \} \cong S^2\ $,
which is a two dimensional round sphere of radius $\ r = 1/\alpha\ $,
we get
\begin{displaymath}
	|x|^2\ \Big| \mean_\mu(0) + 4 \frac{(-x)^\perp}{|x|^2} \Big|
	\geq 2 | (|x - e_3/\alpha| - 1/\alpha , 2 x_4, \ldots, 2 x_n) |
	\geq 2 d(x , S_r)
\end{displaymath}
and again by (\ref{glob.unit-sphere.diam})
and (\ref{glob.unit-sphere.int-fub})
\begin{equation} \label{glob.unit-sphere.dist-sphere}
	\int d(x , S_r )^2 \d \mu(x) \leq \Cn \delta^2.
\end{equation}
Now for $\ x_0 \in spt\ \mu \mbox{ with } d(x_0,M) = 2 \varrho > 0
\mbox{ for } M = span \{ e_1 , e_2 \} \mbox{ or } M = S_r\ $,
we get $\ \mu(B_\varrho(x)) \geq c_0 \varrho^2\ $
by Proposition \ref{haus.meas-low} (\ref{haus.meas-low.esti}),
when observing that $\ spt\ \mu\ $ is connected
and $\ spt\ \mu \not\subseteq B_\varrho(x_0)\ $,
as $\ 0 \not\in B_\varrho(x_0)\ $.
Since obviously $\ d(.,M) \geq \varrho \mbox{ on } B_\varrho(x_0)\ $,
we estimate by (\ref{glob.unit-sphere.dist-plane})
or (\ref{glob.unit-sphere.dist-sphere}) that
\begin{displaymath}
	c_0 \varrho^4
	\leq \int \limits_{B_\varrho(x_0)} d(x,M)^2 \d \mu(x)
	\leq \Cn \delta^2,
\end{displaymath}
hence $\ \varrho \leq \Cn \sqrt{\delta}\ $ and
\begin{displaymath}
	spt\ \mu \subseteq U_{\Cn \sqrt{\delta}}(M)
	:= \{ x \in \rel^n\ |\ d(x,M) < \Cn \sqrt{\delta}\ \}.
\end{displaymath}
For a plane $\ M = span \{ e_1 , e_2 \}\ $,
this is impossible by the next Proposition \ref{glob.sphere}
for $\ \delta\ $ small,
and hence excludes the case with $\ \mean_\mu(0) = 0\ $.
Therefore
\begin{equation} \label{glob.unit-sphere.haus-ink}
	spt\ \mu \subseteq U_{\Cn \sqrt{\delta}}(S_r),
\end{equation}
and we conclude again by the next Proposition \ref{glob.sphere}
and by (\ref{glob.unit-sphere.measure})
for given $\ \varepsilon > 0\ $, if $\ \delta\ $ is small enough, that
\begin{equation} \label{glob.unit-sphere.radius-aux}
%	|\ |\mean_\mu(0)| - 2|,
	|r - 1| \leq \varepsilon.
\end{equation}
After translation to abbreviate notation,
we may assume that $\ S_r
= \partial B_r(0) \cap span \{ e_1 , e_2 , e_3 \}\ $
is centered at the origin.
Then (\ref{glob.unit-sphere.dist-sphere}) reads
\begin{equation} \label{glob.unit-sphere.haus-mod-aux}
	\int |\ |x| - r|^2 \d \mu(x) \leq \Cn \delta^2.
\end{equation}
Returning to the monotinicity in \bcite{kuw.schae.will3} (A.3)
as in (\ref{glob.mini.gamma}) for $\ x = 0\ $
and recalling by \bcite{kuw.schae.will3} (A.7), (A.10), (A.14)
that $\ \gamma(\varrho) \rightarrow \W(\mu)/4
\mbox{ for } \varrho \rightarrow \infty\ $,
we get for large $\ R \mbox{ with } spt\ \mu \subseteq B_R(0)\ $ that
\begin{displaymath}
	\frac{1}{4} \W(\mu)
	= \lim \limits_{\varrho \rightarrow \infty} \gamma(\varrho)
	\geq \gamma(R) =
\end{displaymath}
\begin{displaymath}
	=  R^{-2} \mu(B_R(0))
	+ \frac{1}{16} \int \limits_{B_R(0)}
	|\mean_\mu|^2 \d \mu
	+ \frac{1}{2} \int \limits_{B_R(0)}
	R^{-2}\ x\ \mean_\mu(x) \d \mu(x) =
\end{displaymath}
\begin{displaymath}
	= R^{-2}\ \mu(\rel^n) + \frac{1}{4} \W(\mu)
	+ \frac{1}{2} \int R^{-2}\ x\ \mean_\mu(x) \d \mu(x),
\end{displaymath}
hence by (\ref{glob.unit-sphere.measure}) that
\begin{equation} \label{glob.unit-sphere.hoelder}
	4 \pi
	= \mu(\rel^n)
	\leq - \frac{1}{2} \int x\ \mean_\mu(x) \d \mu(x).
\end{equation}
On the other hand
$\ \parallel \mean_\mu \parallel_{L^2(\mu)} = (4 \W(\mu) )^{1/2}\ $
and by (\ref{glob.unit-sphere.haus-mod-aux})
\begin{equation} \label{glob.unit-sphere.aux}
	\parallel x \parallel_{L^2(\mu)}
	\leq \parallel r \parallel_{L^2(\mu)}
	+ \parallel |x| - r \parallel_{L^2(\mu)}
	\leq r \mu(\rel^n)^{1/2} + \Cn \delta.
\end{equation}
Then we use the H\"older inequality to get
\begin{displaymath}
	1 \geq - \int x\ \mean_\mu(x) \d \mu(x) \Big/
	\Big( \parallel x \parallel_{L^2(\mu)}
	\ \parallel \mean_\mu \parallel_{L^2(\mu)} \Big) \geq
\end{displaymath}
\begin{displaymath}
	\geq \frac{2 \mu(\rel^n)}{(r \mu(\rel^n)^{1/2} + \Cn \delta)
	\ (4 \W(\mu) )^{1/2}}
	\geq \frac{\mu(\rel^n)^{1/2}}{r \mu(\rel^n)^{1/2} + \Cn \delta}
	\geq (1 - \Cn \delta) / r
\end{displaymath}
for $\ \delta\ $ small,
when recalling (\ref{glob.unit-sphere.radius-aux}),
and obtain
\begin{equation} \label{glob.unit-sphere.radius-eta-low}
	r \geq 1 - \Cn \delta.
\end{equation}
To get an estimate from above,
we consider the smooth nearest point projection
$\ \pi: U_{\Cn \sqrt{\delta}}(S_r) \rightarrow S_r\ $
and see that $\ \pi_\# \mu = \theo \Ht \lfloor S_r\ $
is an integral varifold
and $\ \pi_\# T = \thez \curr{S_r}\ $ is an integral current
with measurable $\ \theo: S_r \rightarrow \nat_0,
\thez: S_r \rightarrow \ganz\ $ and 
\begin{equation} \label{glob.unit-sphere.theta}
%\begin{array}{c}
%	\mu_{\pi_\# T} \leq \pi_\# \mu, \\
	\theo = \thez \mbox{ modulo } 2
	\quad \mbox{almost everywhere on } S_r \mbox{ with respect to } \Ht.
%\end{array}
\end{equation}
As $\ \partial \pi_\# T = \pi_\# \partial T = 0\ $,
we see by constancy theorem, see \bcite{sim} Theorem 26.27,
that $\ \thez \in \ganz\ $ is constant.
We claim that
\begin{equation} \label{glob.unit-sphere.theta-unit}
	\theo = 1
	\quad \mbox{on a subset of } S_r \mbox{ with positive measure in } \Ht.
\end{equation}
To this end, we consider $\ \xi \in spt\ \mu \cap S_r\ $
as in (\ref{glob.unit-sphere.int-fub})
and with $\ T_\xi \mu = T_\xi S_r = span \{ e_1 , e_2 \}
\mbox{ and } \mean_\mu(\xi) = 2r^{-1} e_3\ $.
We consider the height-excess
\begin{displaymath}
	heightex_\mu(\xi,\varrho,T_\xi \mu)
	:= \varrho^{-4} \int \limits_{B_\varrho(\xi)}
	dist(x - \xi, T_\xi \mu)^2 \d \mu(x)
	\leq \int \limits_{B_\varrho(\xi)}
	\Big| \frac{(\xi - x)^{\perp_\xi}}{|\xi - x|^2} \Big|^2 \d \mu(x) \leq
\end{displaymath}
\begin{displaymath}
	\leq 2 \int \Big| \frac{1}{4} \mean_\mu(\xi)
	+ \frac{(\xi - x)^{\perp_\xi}}{|\xi-x|^2} \Big|^2
	+ 2 \int \limits_{B_\varrho(\xi)}
	\Big| \frac{1}{4} \mean_\mu(\xi) \Big|^2 \d \mu(x) \leq
\end{displaymath}
\begin{equation} \label{glob.unit-sphere.heightex-aux}
	\leq \delta^2 / 2 + \frac{1}{8} \int \limits_{B_\varrho(\xi)}
	|\mean_\mu(\xi)|^2 \d \mu(x).
\end{equation}
We proceed proving
\begin{equation} \label{glob.unit-sphere.mean-esti}
	\int \limits_A |\mean_\mu|^2 \d \mu
	\leq C \mu(A) + C \delta^2
	\quad \mbox{for any  measurable } A \subseteq spt\ \mu.
\end{equation}
First we assume that $\ 2 diam\ A \leq diam\ spt\ \mu =: d\ $.
Then there exists $\ x \in spt\ \mu \mbox{ with } d(x,A) \geq d/2\ $,
and we get from (\ref{glob.unit-sphere.diam}) and (\ref{glob.unit-sphere.int})
\begin{displaymath}
	\int \limits_A |\mean_\mu|^2 \d \mu
	\leq 32 \int \limits_A
	\Big| \frac{(y - x)^{\perp_y}}{|y-x|^2} \Big|^2 \d \mu(y)
	+ 8 \delta^2
	\leq 128 d^{-2} \mu(A) + 8 \delta^2
	\leq C \mu(A) + 8 \delta^2.
\end{displaymath}
In the general case
we take a maximal subset $\ \{ x_1, \ldots , x_N \} \subseteq spt\ \mu
\mbox{ with } |x_i -x_j| \geq d/4 \mbox{ for } i \neq j\ $.
Then $\ B_{d/8}(x_i)\ $ are pairwise disjoint
and $\ spt\ \mu \not\subseteq B_{d/8}(x_i)\ $,
and we get by (\ref{glob.unit-sphere.measure}),
(\ref{haus.meas-low.esti})
and connectedness of $\ spt\ \mu\ $ that
\begin{displaymath}
	4 \pi = \mu(\rel^n) \geq \sum \limits_{i=1}^N \mu(B_{d/8}(x_i))
	\geq c_0 N d^2 / 16,
\end{displaymath}
hence $\ N \leq C\ $ again by (\ref{glob.unit-sphere.diam}).
Since on the other hand $\ spt\ \mu \subseteq \cup_{i=1}^N B_{d/4}(x_i)\ $,
we get from above
\begin{displaymath}
	\int \limits_A |\mean_\mu|^2 \d \mu
	\leq \sum \limits_{i=1}^N
	\int \limits_{A \cap B_{d/4}(x_i)}|\mean_\mu|^2 \d \mu
	\leq \sum \limits_{i=1}^N \Big( C \mu(A) + 8 \delta^2 \Big)
	\leq C \mu(A) + C \delta^2,
\end{displaymath}
which is (\ref{glob.unit-sphere.mean-esti}).

Then we obtain from (\ref{glob.unit-sphere.heightex-aux})
combined with (\ref{glob.unit-sphere.density})
\begin{displaymath}
	heightex_\mu(\xi,\varrho,T_\xi \mu) \leq
\end{displaymath}
\begin{equation} \label{glob.unit-sphere.heightex}
	\leq \delta^2 / 2 + \frac{1}{8} \int \limits_{B_\varrho(\xi)}
	|\mean_\mu(\xi)|^2 \d \mu(x)
	\leq C \mu(B_\varrho(\xi)) + C \delta^2
	\leq C (\varrho^2 + \delta^2)
\end{equation}
and for the tilt-excess by \bcite{sim} Lemma 22.2 that
\begin{displaymath}
	tiltex_\mu(\xi,\varrho,T_\xi \mu)
	:= \varrho^{-2} \int \limits_{B_\varrho(\xi)}
	\parallel T_x \mu , T_\xi \mu \parallel^2 \d \mu(x) \leq
\end{displaymath}
\begin{equation} \label{glob.unit-sphere.tiltex}
	\leq C \Big( heightex_\mu (\xi,2 \varrho,T_\xi \mu)
	+ \parallel \mean_\mu \parallel_{L^2(\mu,B_{2 \varrho}(\xi))}^2 \Big)
%	heightex_\mu (\xi,2 \varrho,T_\xi \mu)^{1/2} \Big)
	\leq C (\varrho^2 + \delta^2).
\end{equation}
Next (\ref{glob.unit-sphere.mean-esti}) 
yields by \bcite{kuw.schae.will3} (A.6) and (A.10)
for any $\ 0 < \tau < 1/2\ $ that
\begin{displaymath}
	\frac{\mu(B_\varrho(\xi))}{\pi \varrho^2}
	\geq (1 + \tau)^{-1} \theta^2(\mu,\xi)
	- C (1 + \tau^{-1}) \int \limits_{B_\varrho(\xi)}
	|\mean_\mu(\xi)|^2 \d \mu(x) \geq
\end{displaymath}
\begin{displaymath}
	\geq 1 - \tau - C \tau^{-1} (\varrho^2 + \delta^2),
\end{displaymath}
hence
\begin{equation} \label{glob.unit-sphere.density-low}
	\frac{\mu(B_\varrho(\xi))}{\pi \varrho^2}
	\geq 1 - C \varrho
	\quad \mbox{for } \varrho \geq \delta.
\end{equation}
Moreover by \bcite{kuw.schae.will3} (A.4), (A.5),
the monotonicity of $\ \gamma\ $
that $\ \gamma(\varrho) \rightarrow \W(\mu)/4
\mbox{ for } \varrho \rightarrow \infty\ $,
by \bcite{kuw.schae.will3} (A.7), (A.10), (A.14)
\begin{displaymath}
	\frac{\mu(B_\varrho(\xi))}{\pi \varrho^2}
	\leq W(\mu) / (4 \pi) + |R_{\xi,\varrho}| / \pi \leq
\end{displaymath}
\begin{displaymath}
	\leq 1 + \delta^2 / (4 \pi) + (1/(2 \pi))
	(\varrho^{-2} \mu(B_\varrho(\xi)))^{1/2}
	\ \parallel \mean_\mu \parallel_{L^2(\mu,B_{2 \varrho}(\xi))}
\end{displaymath}
\begin{displaymath}
	\leq 1 + \delta^2 / (4 \pi)
	+ \tau \frac{\mu(B_\varrho(\xi))}{\pi \varrho^2}
	+ C \tau^{-1}  \parallel \mean_\mu
	\parallel_{L^2(\mu,B_{2 \varrho}(\xi))}^2 \leq
\end{displaymath}
\begin{displaymath}
	\leq 1 + C \tau^{-1} (\varrho^2 + \delta^2)
	+ \tau \frac{\mu(B_\varrho(\xi))}{\pi \varrho^2}
\end{displaymath}
for any $\ 0 < \tau < 1\ $, hence
\begin{equation} \label{glob.unit.sphere.density-upp}
	\frac{\mu(B_\varrho(\xi))}{\pi \varrho^2}
	\leq 1 + C \varrho
	\quad \mbox{for } \varrho \geq \delta.
\end{equation}
Combining (\ref{glob.unit-sphere.mean-esti}),
(\ref{glob.unit-sphere.heightex})
and (\ref{glob.unit-sphere.tiltex}),
we get from \bcite{brakke} Theorem 5.4
for $\ \varrho \geq \delta\ $ small enough
a lipschitz approximation
of $\ \mu \mbox{ over } T_\xi \mu
\mbox{ at } \xi = (\xi',\xi'') \in \rel^2 \times \rel^{n-2}\ $,
that is there exists a single-valued lipschitz map
\begin{displaymath}
\begin{array}{c}
	f: B^2_\varrho(\xi') \subseteq T_\xi \mu = span \{ e_1, e_2 \}
	\rightarrow T^\perp_\xi \mu = span \{ e_3, \ldots , e_n \}, \\
	F: B^2_\varrho(\xi') \rightarrow \rel^2 \times \rel^{n-2},
	\quad F(y) = (y,f(y)), \\
\end{array}
\end{displaymath}
satisfying
\begin{equation} \label{glob.unit-sphere.lip.lip}
\begin{array}{c}
	lip\ f \leq 1, \\
	\varrho^{-1} \parallel f - \xi''
	\parallel_{L^\infty(B^2_\varrho(\xi'))}
	\leq \Cn \varrho^{2/(n+2)} \leq 1, \\
\end{array}
\end{equation}
for $\ \varrho\ $ small enough,
and there exists a Borel set $\ Y \subseteq B^2_\varrho(\xi')\ $ such that
\begin{equation} \label{glob.unit-sphere.lip.except.y}
	\theta^2(\mu,(y,z)) = \delta_{f(y) z}
	\quad \mbox{for all } y \in Y \subseteq \rel^2,
	\ z \in B^{n-2}_{1}(\xi'') \subseteq \rel^{n-2},
\end{equation}
and setting
\begin{equation} \label{glob.unit-sphere.lip.except.x}
	X := spt\ \mu \cap (Y \times B^{n-2}_\varrho(\xi'')) = F(Y),
\end{equation}
then
\begin{equation} \label{glob.unit-sphere.lip.esti}
	\varrho^{-2} \mu((B^2_\varrho(\xi')
	\times B^{n-2}_\varrho(\xi'')) - X) +
	\varrho^{-2} \Lt(B^2_\varrho(\xi') - Y) \leq \Cn \varrho^2.
\end{equation}
Now for the nearest point projection
$\ |\pi(x) - x| \leq d(x,S_r) \leq \Cn \sqrt{\delta}\ $
by (\ref{glob.unit-sphere.haus-ink}),
hence for $\ \varrho \geq C_n \sqrt{\delta}\ $ that
\begin{equation} \label{glob.unit-sphere.lip.aux}
\begin{array}{c}
	\pi(spt\ \mu - B_\varrho(\xi)) \cap B_{\varrho/2}(\xi) = \emptyset, \\
	\pi(spt\ \mu \cap B_{\varrho/4}(\xi)) \subseteq B_{\varrho/2}(\xi), \\
\end{array}
\end{equation}
Next $\ \pi \mbox{ is injective on } X = F(Y)\ $,
indeed for $\ y_1,y_2 \in Y
\mbox{ with } \pi(F(y_1)) = \pi(F(y_2)) =: p \in S_r\ $,
we have $\ F(y_2) - F(y_1) \in N_p S_r\ $,
which denotes the normal space.
As $\ p \in B_{2 \varrho + \Cn \sqrt{\delta}}(\xi)
\mbox{ and } T_\xi S_r = T_\xi \mu = span \{ e_1 , e_2 \}\ $, we get
\begin{displaymath}
	|y_2 -y_1| = | \pi_{\rel^2 \times \{0\} } (F(y_2) - F(y_1)) |
%	= | \pi_{T_\xi S_r} \pi_{N_p S_r} (F(y_2) - F(y_1)) |
	=
\end{displaymath}
\begin{displaymath}
	= | (\pi_{T_\xi S_r} - \pi_{T_p S_r}) (F(y_2) - F(y_1)) |
	\leq \frac{1}{4} |F(y_2) - F(y_1)|
	\leq \frac{1}{2} |y_2 - y_1|
\end{displaymath}
for $\ \varrho\ $ small enough
and as $\ lip\ f \leq 1\ $ by (\ref{glob.unit-sphere.lip.lip}).
This implies $\ y_1 = y_2\ $,
and $\ \pi \mbox{ is injective on } X\ $.
Since further $\ \theta^2(\mu) = 1 \mbox{ on } X\ $
by (\ref{glob.unit-sphere.lip.except.y}),
we get $\ \pi_\# (\mu \lfloor X) = \Ht \lfloor \pi(X)\ $,
hence by (\ref{glob.unit-sphere.lip.aux})
\begin{equation} \label{glob.unit-sphere.lip.aux2}
	\theo = \theta^2(\pi_\# \mu) = 1
	\mbox{ on } \pi(X) \cap B_{\varrho/2}(\xi)
	- \pi(spt\ \mu \cap B_\varrho(\xi) - X).
\end{equation}
Clearly $\ \pi\ $ is lipschitz with some uniform constant $\ L < \infty\ $,
when observing $\ r \approx 1\ $ by (\ref{glob.unit-sphere.radius-aux})
and for $\ \delta\ $ small, hence we have by (\ref{glob.unit-sphere.lip.esti})
when observing that $\ \mu = \Ht \lfloor spt\ \mu\ $
by (\ref{glob.unit-sphere.support}) that
\begin{displaymath}
	\varrho^{-2} \Ht( \pi(spt\ \mu \cap B_\varrho(\xi) - X) )
	\leq L^2 \varrho^{-2} \Ht( spt\ \mu \cap B_\varrho(\xi) - X ) =
\end{displaymath}
\begin{displaymath}
	= L^2 \varrho^{-2} \mu(B_\varrho(\xi) - X)
	\leq \Cn \varrho^2.
\end{displaymath}
On the other hand $\ (J_\mu \pi)= J_{T \mu} \pi \geq c_0 > 0
\mbox{ on } X\ $ for $\ \varrho\ $ small,
as $\ lip\ f \leq 1\ $ by (\ref{glob.unit-sphere.lip.lip}),
and recalling $\ \pi_\# \mu = \pi(J_\mu \pi \cdot \mu)\ $
we get by (\ref{glob.unit-sphere.density-low}),
(\ref{glob.unit-sphere.lip.esti})
and (\ref{glob.unit-sphere.lip.aux})
\begin{displaymath}
	\Ht(\pi(X) \cap B_{\varrho/2}(\xi))
	= (\pi_\# \mu)(\pi(X) \cap B_{\varrho/2}(\xi))
	= \pi(J_\mu \pi \cdot \mu)(\pi(X) \cap B_{\varrho/2}(\xi)) =
\end{displaymath}
\begin{displaymath}
	= (J_\mu \pi \cdot\mu)(X \cap \pi^{-1}(B_{\varrho/2}(\xi)))
	\geq c_0 \mu(X \cap B_{\varrho/4}(\xi))
	\geq c_0 ( \mu(B_{\varrho/4}(\xi))
	- \mu( B_{\varrho/4}(\xi) - X) ) \geq
\end{displaymath}
\begin{displaymath}
	\geq c_0 \varrho^2 ( (1/16) - C \varrho - \Cn \varrho^2).
\end{displaymath}
Together we see
\begin{displaymath}
	\Ht \Big( \pi(X) \cap B_{\varrho/2}(\xi)
	- \pi(spt\ \mu \cap B_\varrho(\xi) - X) \Big)
	\geq c_0 \varrho^2
	( (1/16) - C \varrho - C \varrho^2 - \Cn \varrho^2) > 0
\end{displaymath}
for $\ \varrho \geq \Cn \sqrt{\delta}\ $ small enough,
which yields (\ref{glob.unit-sphere.theta-unit}).

Then we conclude by (\ref{glob.unit-sphere.theta}) that $\ \thez \neq 0\ $,
hence $\ \Ht(S_r - \pi(spt\ \mu)) = 0\ $, in particular
\begin{equation} \label{glob.unit-sphere.surjective}
	\Ht(\pi(spt\ \mu)) = \Ht(S_r) = 4 \pi r^2
\end{equation}
and by (\ref{glob.unit-sphere.haus-ink})
\begin{equation} \label{glob.unit-sphere.haus-sup}
	S_r \subseteq U_{\Cn \sqrt{\delta}}(spt\ \mu).
\end{equation}
On the other hand, we have for $\ x = (u,v) \in U_{\Cn \sqrt{\delta}}(S_r)
\subseteq \rel^3 \times \rel^{n-3}\ $ that $\ \pi(u,v) = r u / |u|\ $ and
\begin{displaymath}
	(D \pi)(u,v) = (r/|u|) \Big( I_3 - u^T u / |u|^2 , 0 \Big),
\end{displaymath}
hence for the Jacobian
\begin{displaymath}
	J_\mu \pi(x) \leq \parallel D \pi(x) \parallel^2 \leq r^2 / |u|^2.
\end{displaymath}
Then by the Area formula, see \bcite{sim} (12.4), and (\ref{glob.unit-sphere.support})
\begin{displaymath}
	\Ht(\pi(spt\ \mu)) \leq \int \limits_{spt\ \mu} J_\mu \pi \d \Ht
	\leq \int r^2 / |u|^2 \d \mu(x)
\end{displaymath}
and combing with (\ref{glob.unit-sphere.measure})
and (\ref{glob.unit-sphere.surjective})
\begin{equation} \label{glob.unit-sphere.area}
	\int 1 \d \mu(x) = 4 \pi \leq \int 1 / |u|^2 \d \mu(x).
\end{equation}
Estimating
\begin{equation} \label{glob.unit-sphere.jac-aux}
	| r - |u| | \leq |\pi(x) - u| \leq |\pi(x) - x| = d(x,S_r)
\end{equation}
and
\begin{displaymath}
	1 / |u|^2 \leq 1 / (r - d(x,S_r))^2 \leq r^{-2} (1 + C d(x,S_r))
\end{displaymath}
for $\ \delta \mbox{ and } \varepsilon\ $
small enough and taking into account that
$\ x \in U_{\Cn \sqrt{\delta}}(S_r) \ $ and (\ref{glob.unit-sphere.radius-aux}),
we get from (\ref{glob.unit-sphere.area})
when using (\ref{glob.unit-sphere.dist-sphere})
\begin{displaymath}
	4 \pi r^2
	\leq \int (1 + C d(x,S_r)) \d \mu(x)
%	\leq
%\end{displaymath}
%\begin{displaymath}
%	\leq \mu(\rel^n) + C \Big( \int d(x,S_r)^2 \d \mu(x) \Big)^{1/2}
	\leq 4 \pi + \Cn \delta,
\end{displaymath}
hence $\ r^2 \leq 1 + \Cn \delta
\mbox{ and } r \leq 1 + \Cn \delta\ $
and with (\ref{glob.unit-sphere.radius-eta-low})
\begin{equation} \label{glob.unit-sphere.radius-eta}
	|r - 1| \leq \Cn \delta.
\end{equation}
Putting $\ S_1 := \partial B_1(0) \cap span \{ e_1 , e_2 , e_3 \}\ $,
we obtain (\ref{glob.unit-sphere.sphere})
from (\ref{glob.unit-sphere.dist-sphere})
and (\ref{glob.unit-sphere.radius-eta}).
Combining (\ref{glob.unit-sphere.haus-ink}),
(\ref{glob.unit-sphere.haus-sup})
and (\ref{glob.unit-sphere.radius-eta})
yields (\ref{glob.unit-sphere.haus}).
%Clearly (\ref{glob.unit-sphere.haus-mod-aux})
%and (\ref{glob.unit-sphere.radius-eta}) give
%\begin{equation} \label{glob.unit-sphere.haus-mod}
%	\int |\ |x| - 1|^2 \d \mu(x) \leq C \delta^2.
%\end{equation}

Finally to prove (\ref{glob.unit-sphere.mean}),
we return to (\ref{glob.unit-sphere.jac-aux})
and see using (\ref{glob.unit-sphere.radius-eta})
\begin{displaymath}
	|1 - |u| | \leq |r - |u| | + \Cn \delta
	\leq d(x,S_r) + \Cn \delta
	\leq \Cn \sqrt{\delta} \leq 1/2
\end{displaymath}
for $\ \delta\ $ small, hence
\begin{displaymath}
	1 / |u|^2 = (1 + |u|^2 - 1)^{-1}
	\leq 1 - (|u|^2 - 1) + C (|u|^2 - 1)^2 \leq
\end{displaymath}
\begin{displaymath}
	\leq 1 - (|x|^2 - 1) + |v|^2 + C (|u| - 1)^2
	\leq 1 - (|x|^2 - 1) + C d(x,S_r)^2 + \Cn \delta^2
\end{displaymath}
when observing that $\ |v| \leq d(x,S_r)\ $.
Plugging into (\ref{glob.unit-sphere.area}), we get
\begin{displaymath}
	\int 1 \d \mu(x)
	\leq \int \Big( 1 - (|x|^2 - 1)
	+ C d(x,S_r)^2 + \Cn \delta^2 \Big) \d \mu(x)	
\end{displaymath}
and using (\ref{glob.unit-sphere.dist-sphere})
\begin{displaymath}
	\int |x|^2 \d \mu(x) \leq 4 \pi + \Cn \delta^2.
\end{displaymath}
Then by (\ref{glob.unit-sphere.energy}) and (\ref{glob.unit-sphere.hoelder})
\begin{displaymath}
	\int \Big| \frac{1}{2} \mean_\mu(x) + x \Big|^2 \d \mu(x)
	= \W(\mu)
%	= \frac{1}{4} \parallel \mean_\mu \parallel_{L^2(\mu)}^2
	+ \int \mean_\mu(x)\ x \d \mu(x)
	+ \int |x|^2 \mu(x) \leq
\end{displaymath}
\begin{displaymath}
	\leq 4 \pi + \delta^2  - 8 \pi + 4 \pi + \Cn \delta^2 \leq \Cn \delta^2,
\end{displaymath}
which yields (\ref{glob.unit-sphere.mean}).
\proof

\begin{proposition} \label{glob.sphere}

Let $\ \mu \neq 0\ $ be an integral $\ 2-$varifold
with square integrable weak mean curvature,
compact support,
\begin{equation} \label{glob.sphere.measure}
	\mu(\rel^n) = 4 \pi
\end{equation}
and
\begin{equation} \label{glob.sphere.energy}
	\W(\mu) \leq 8 \pi - \tau
\end{equation}
for some $\ \tau > 0\ $.
Then for any $\ \varepsilon > 0\ $
there exists $\ \eta = \eta(\tau,\varepsilon)\ $
such that if $\ spt\ \mu\ $ is close to a two dimensional round sphere
or two dimensional plane, more precisely after rotation
\begin{equation} \label{glob.sphere.sphere}
	spt\ \mu \subseteq U_\eta \Big( \partial B_r(0)
	\cap (\rel^3 \times \{0\}) \Big)
\end{equation}
or
\begin{equation} \label{glob.sphere.plane}
	spt\ \mu \subseteq U_\eta(\rel^2 \times \{0\}),
\end{equation}
which corresponds to $\ r = \infty\ $
in the sense that (\ref{glob.sphere.plane})
implies (\ref{glob.sphere.sphere}) for some large $\ r\ $,
as $\ spt\ \mu\ $ is compact,
then
\begin{equation} \label{glob.sphere.radius}
	|r - 1| \leq \varepsilon,
\end{equation}
in particular the second case is excluded,
and there exists a two dimensional round sphere
$\ S_1 \subseteq \rel^n\ $ with radius 1 and
\begin{equation} \label{glob.sphere.haus}
	d_H(spt\ \mu , S_1) \leq \varepsilon,
\end{equation}
where $\ d_H\ $ denotes the Hausdorff distance.
\end{proposition}
{\pr Proof:} \\
Firstly (\ref{glob.sphere.energy}) implies by (\ref{haus.mono.conn})
that $\ spt\ \mu\ $ is connected.
Next by (\ref{glob.sphere.measure}),
we get from \bcite{sim.will} Lemma 1.1
and (\ref{haus.meas-low.esti}) by connectedness of $\ spt\ \mu\ $
for the non-smooth case that
\begin{equation} \label{glob.sphere.diam}
	diam\ spt\ \mu \leq \Cn
\end{equation}
for some $\ \Cn < \infty\ $.
Therefore in case of (\ref{glob.sphere.plane}),
we get (\ref{glob.sphere.sphere}) for $\ \eta \mbox{ replaced by } 2 \eta\ $
and some $\ r \gg 1\ $.

Suppose (\ref{glob.sphere.radius}) is wrong for some $\ \varepsilon > 0\ $,
then there exists a sequence of varifolds $\ \mu_m\ $ as above
with $\ \W(\mu_m) \leq 8 \pi - \tau\ $
and radii $\ r_m \mbox{ with } |r_m - 1| \geq \varepsilon\ $
and after approprite rotation
\begin{equation} \label{glob.sphere.close}
	spt\ \mu_m \subseteq U_{1/m} \Big( \partial B_{r_m}(0)
	\cap span \{ e_1 , e_2 , e_3 \} \Big),
\end{equation}
Firstly by \bcite{kuw.schae.will3} (A.16)
\begin{equation} \label{glob.sphere.density}
	\varrho^{-2} \mu_m(B_\varrho) \leq C
	\quad \mbox{ for all } B_\varrho \subseteq \rel^n,
\end{equation}
hence by (\ref{glob.sphere.measure})
\begin{equation} \label{glob.sphere.diam-low}
	diam\ spt\ \mu_m \geq c_0
\end{equation}
for some $\ c_0 > 0\ $, in particular by (\ref{glob.sphere.close})
\begin{equation} \label{glob.sphere.radius-low}
	\liminf \limits_{m \rightarrow \infty} r_m > 0.
\end{equation}
Next if $\ r_m\ $ were unbounded,
we get $\ r_m \rightarrow \infty\ $ after passing to a subsequence.
We select $\ x_m \in spt\ \mu_m \neq \emptyset\ $
and assume by (\ref{glob.sphere.close}) after a further rotation
that $\ |x_m - r_m e_3| < 1/m\ $.
Translating to $\ \nu_m := (x \mapsto x - r_m e_3)_\# \mu_m\ $,
we clearly have
\begin{displaymath}
\begin{array}{c}
	0 \leftarrow y_m := x_m - r_m e_3 \in spt\ \nu_m, \\
	spt\ \nu_m \subseteq U_{1/m} \Big( \partial B_{r_m}(r_m e_3)
	\cap span \{ e_1 , e_2 , e_3 \} \Big),
\end{array}
\end{displaymath}
by (\ref{glob.sphere.close})
and get from (\ref{glob.sphere.density}) after passing to a subsequence
that $\ \nu_m \rightarrow \nu\ $ weakly as Radon measures.
As the supports of $\ \mu_m \mbox{ respectively of } \nu_m\ $
are connected and their diameters
are bounded from below by (\ref{glob.sphere.diam-low}),
we get from Proposition \ref{haus.prop}
that $\ spt\ \nu_m \rightarrow spt\ \nu\ $ locally in Hausdorff distance,
hence
\begin{displaymath}
	0 \in spt\ \nu \subseteq span \{ e_1 , e_2 \}.
\end{displaymath}
Next by boundedness $\ \W(\nu_m) = \W(\mu_m) \leq 8 \pi\ $,
we get that $\ \nu\ $ has weak mean curvature in $\ L^2(\nu)\ $
and by Allard's integral compactness theorem,
see \bcite{allard.vari} Theorem 6.4 or \bcite{sim} Remark 42.8,
that $\ \nu\ $ is an integral varifold.
Clearly $\ T \nu = span \{ e_1 , e_2 \}\ $
and further $\ \mean_\nu \perp T \nu\ $
by \bcite{brakke} Theorem 5.8
almost everywhere with respect to $\ \nu\ $.
Then $\ \nu\ $ is stationary in $\ span \{ e_1 , e_2 \}\ $
in the sense of \bcite{sim} 16.4 or 41.2 (3),
and by constancy theorem, see \bcite{sim} Theorem 41.1,
we get $\ \nu = \theta \Ht \lfloor span \{ e_1 , e_2 \}
\mbox{ for some constant } \theta > 0\ $,
recalling that $\ \nu \neq 0\ $ by above.
Therefore $\ \nu(\rel^n) = \infty\ $,
but $\ \nu(\rel^n)
\leq \limsup_{m \rightarrow \infty} \nu_m(\rel^n)
= 4 \pi\ $ by (\ref{glob.sphere.measure}),
which is a contradiction,
hence $\ r_m\ $ is bounded.

Combining with (\ref{glob.sphere.radius-low}),
we get after passing to a subsequence
\begin{equation} \label{glob.sphere.radius-conv}
	r_m \rightarrow r > 0, r \neq 1,
\end{equation}
when recalling that $\ |r_m-1| \geq \varepsilon\ $ by assumption.

Then as above we get after passing to a further subsequence
that $\ \mu_m \rightarrow \mu\ $ weakly as Radon measures
with $\ \mu\ $ an integral varifold
with square integrable weak mean curvature and
\begin{displaymath}
	spt\ \mu \subseteq \partial B_r(0) \cap span \{ e_1 , e_2 , e_3 \} =: S_r \cong S^2,
\end{displaymath}
hence with constancy theorem, see \bcite{sim} Theorem 41.1,
we get $\ \mu = \theta \Ht \lfloor S_r
\mbox{ for some constant } \theta \in \nat_0\ $,
as $\ \mu\ $ is integral.
We conclude
\begin{equation} \label{glob.sphere.meas-a}
	\mu(\rel^n) = 4 \pi \theta r^2.
\end{equation}
Moreover as $\ spt\ \mu \subseteq B_{r_m}(0) \subseteq B_{r+1}(0)
\mbox{ for large } m\ $, we get  by (\ref{glob.sphere.measure})
\begin{equation} \label{glob.sphere.meas-b}
	\mu(\rel^n) \leftarrow \mu_m(\rel^n) = 4 \pi
\end{equation}
and by lower semicontinuity and (\ref{glob.sphere.energy})
\begin{displaymath}
	\W(\mu) \leq \limsup \limits_{m \rightarrow \infty} \W(\mu_m) < 8 \pi.
\end{displaymath}
Therefore $\ \mu \neq 0\ $, hence $\ \theta > 0\ $,
and secondly by the Li-Yau inequality,
see \bcite{kuw.schae.will3} (A.17),
we get $\ \theta^2(\mu) \leq \W(\mu) / 4 \pi < 2\ $,
hence $\ \theta < 2\ $.
As $\ \theta\ $ is an integer, we conclude $\ \theta = 1\ $.
Combining (\ref{glob.sphere.meas-a}) and (\ref{glob.sphere.meas-b}),
we see $\ r = 1\ $,
which contradicts (\ref{glob.sphere.radius-conv}),
and (\ref{glob.sphere.radius}) is proved.

The same procedure yields (\ref{glob.sphere.haus}).
Indeed we see that $\ r = 1 \mbox{ and } S_1\ $
is a two dimensional round sphere of radius 1.
Moreover $\ \mu_m \rightarrow \mu = \Ht \lfloor S_1\ $
and by Proposition \ref{haus.prop} the convergence
$\ spt\ \mu_m \rightarrow spt\ \mu = S_1\ $ is in global Hausdorff distance,
as $\ spt\ \mu_m\ $ stay inside
a fixed bounded domain by (\ref{glob.sphere.diam}).
Then for large $\ m\ $, we get (\ref{glob.sphere.haus}) with $\ S_1\ $,
which finishes the proof of the proposition.
\proof
%%%%%
Up to this point all of our arguments worked for integral $2$-varifolds with square integrable weak mean curvature. From now on we work with smooth surfaces,
even though one can relax this regularity assumption slightly,
see the remark after Theorem \ref{intro.theorem-funda}.

\begin{proposition} \label{glob-smooth.sphere}

Let $\ \Sigma \subseteq \rel^n\ $
be a smooth embedded closed surface with
\begin{equation} \label{glob-smooth.sphere.area}
	\Ht(\Sigma) = 4 \pi
\end{equation}
and
\begin{equation} \label{glob-smooth.sphere.funda-0}
	\parallel A^0_\Sigma \parallel_{L^2(\Sigma)} \leq \delta.
\end{equation}
Then for $\ \delta\ $ small depending on $\ n\ $
\begin{equation} \label{glob-smooth.sphere.funda}
	\parallel A_\Sigma - \nor g \parallel_{L^2(\Sigma)}
	\leq \Cn \delta
\end{equation}
for some measurable unit normal vector field
$\ \nor \mbox{ on } \Sigma\ $ and
\begin{equation} \label{glob-smooth.sphere.gauss}
	\parallel K_\Sigma - 1 \parallel_{L^1(\Sigma)}
	\leq \Cn \delta.
\end{equation}
Moreover after an appropriate translation and rotation
\begin{equation} \label{glob-smooth.sphere.mean}
	\parallel \mean_\Sigma + 2 id_\Sigma \parallel_{L^2(\Sigma)}
	\leq \Cn \delta,
\end{equation}
\begin{equation} \label{glob-smooth.sphere.sphere}
	\parallel |id_\Sigma| - 1 \parallel_{L^2(\Sigma)}
	\leq \Cn \delta,
\end{equation}
\begin{equation} \label{glob-smooth.sphere.l-infty}
	\parallel id_\Sigma \parallel_{L^\infty(\Sigma)} \leq 2.
\end{equation}
\end{proposition}
{\pr Proof:} \\
First we see for $\ \delta^2 < 4 \pi\ $
by the Gau\ss\ equations
and the Gau\ss-Bonnet theorem in (\ref{intro.gauss})
that $\ \chi(\Sigma) = 2\ $ and $\ \Sigma\ $ is a sphere,
in particular $\ \Sigma\ $ is orientable.
We put $\ \mu = \Ht \lfloor \Sigma, T := \mu \llcorner \xi\ $
for some smooth orientation $\ \xi \mbox{ on } \Sigma\ $
and see
\begin{displaymath}
\begin{array}{c}
	\mu(\rel^n) = \Ht(\Sigma) = 4 \pi,
	\quad \mean_\mu = \mean_\Sigma, \\
	\mu_T = \mu, \quad \partial T = 0, \\
\end{array}
\end{displaymath}
and by the Gau\ss\ equations and the Gau\ss-Bonnet theorem
in (\ref{intro.gauss}) that
\begin{equation} \label{glob-smooth.sphere.will}
	\W(\mu) = \W(\Sigma)
%	= \frac{1}{2} \int \limits_{\Sigma}
%	|A^0_\Sigma|^2 \d \Ht + 2 \pi \chi(\Sigma)
	\leq 4 \pi + \delta^2 / 2.
\end{equation}
since $\ \chi(\Sigma) \leq 2\ $.
Then
after an appropriate translation and rotation
Proposition \ref{glob.unit-sphere} implies
(\ref{glob-smooth.sphere.mean}),
(\ref{glob-smooth.sphere.sphere})
and (\ref{glob-smooth.sphere.l-infty}).
This yields
\begin{equation} \label{glob-smooth.sphere.mean-aux}
	\Big\|\ \Big| \frac{1}{2} \mean \Big| - 1 \Big\|_{L^2(\Sigma)}
	\leq \Big\|\ \frac{1}{2} \mean + id_\Sigma \Big\|_{L^2(\Sigma)}
	+ \parallel |id_\Sigma| - 1 \parallel_{L^2(\Sigma)}
	\leq \Cn \delta,
\end{equation}
and putting $\ \nor := \mean / |\mean| \mbox{ for } \mean \neq 0
\mbox{ and } \nor \in N \Sigma, |\nor| = 1
\mbox{ and } \nor\ $ measurable otherwise,
we see
\begin{displaymath}
	\parallel \frac{1}{2} \mean - \nor \parallel_{L^2(\Sigma)}
	= \Big\|\ \Big| \frac{1}{2} \mean \Big| - 1 \Big\|_{L^2(\Sigma)}
	\leq \Cn \delta,
\end{displaymath}
hence
\begin{displaymath}
	\parallel A - \nor g \parallel_{L^2(\Sigma)}
	\leq \parallel A^0 \parallel_{L^2(\Sigma)}
	+ \sqrt{2} \parallel \frac{1}{2} \mean - \nor \parallel_{L^2(\Sigma)}
	\leq \Cn \delta,
\end{displaymath}
which is (\ref{glob-smooth.sphere.funda}).

Next by the general estimate
of the Gau\ss-curvature in the next proposition,
we get combining with (\ref{glob-smooth.sphere.area}),
(\ref{glob-smooth.sphere.will})
and (\ref{glob-smooth.sphere.mean-aux})
\begin{displaymath}
	\parallel K_\Sigma - 1 \parallel_{L^1(\Sigma)}
	\leq \Big\| K - \Big|\frac{1}{2} \mean \Big|^2 \Big\|_{L^1(\Sigma)}
	+ \Big\|\ \Big| \frac{1}{2} \mean \Big|^2 - 1 \Big\|_{L^1(\Sigma)} \leq
\end{displaymath}
\begin{displaymath}
	\leq \Cn \parallel A^0_\Sigma \parallel_{L^2(\Sigma)}^2
	+ \Big\|\ \Big| \frac{1}{2} \mean \Big| + 1 \Big\|_{L^2(\Sigma)}
	\ \Big\|\ \Big| \frac{1}{2} \mean \Big| - 1 \Big\|_{L^2(\Sigma)}
	\leq \Cn \delta.
\end{displaymath}
which is (\ref{glob-smooth.sphere.gauss}).
\proof
{\large \bf Remark:} \\
Assuming embeddedness in the above proposition is no restriction,
as we see from (\ref{glob-smooth.sphere.will})
that $\ \W(\Sigma) < 8 \pi \mbox{ for } \delta^2 < 8 \pi\ $,
and embeddedness follows from the Li-Yau inequality in \bcite{li.yau}
or \bcite{kuw.schae.will3} (A.17).
\defin
\\ \ \\
We supplement the general estimate of the Gau\ss-curvature.

\begin{proposition} \label{glob-smooth.gauss}

Let $\ f: \Sigma \rightarrow \rel^n\ $ be a smooth immersion
of an open surface $\ \Sigma\ $.
Then the Gau\ss-curvature $\ K \mbox{ of } f\ $ is estimated by
\begin{displaymath}
	\Big| K - \Big|\frac{1}{2} \mean \Big|^2 \Big|
	\leq \Cn |A^0|^2_g.
\end{displaymath}
\end{proposition}
{\pr Proof:} \\
We fix $\ p \in \Sigma\ $
and assume in local charts that $\ g_{ij}(p) = \delta_{ij}\ $
and choose an orthonormal basis $\ v_1, \ldots, v_{n-2}\ $
of the normal space $\ N_p f \mbox{ of } f \mbox{ at } p\ $
with $\ \mean = |\mean| v_1 =: 2 \alpha v_1\ $.
The Gau\ss\ curvature
can be written by the Gau\ss\ equations,
see \bcite{docarmo} \S 6 Proposition 3.1, as
\begin{displaymath}
        K = \langle A_{11} , A_{22} \rangle
        - \langle A_{12} , A_{21} \rangle
	= \sum \limits_{m=1}^{n-2} \det(A_{ij} \cdot v_m).
\end{displaymath}
Clearly $\ A_{ij} \cdot v_m = A^0_{ij} \cdot v_m\ $,
hence for $\ m \geq 2\ $
\begin{displaymath}
	\Big| \sum \limits_{m=2}^{n-2} \det(A_{ij} \cdot v_m) \Big|
	\leq \Cn |A^0|^2_g.
\end{displaymath}
Next
\begin{displaymath}
	A_{ij} \cdot v_1
	= \frac{1}{2} \mean \cdot v_1\ \delta_{ij} + A^0_{ij} \cdot v_1
	= \alpha \delta_{ij} + A^0_{ij} \cdot v_1
\end{displaymath}
and abbreviating $\ h^{0,1}_{ij} := A^0_{ij} \cdot v_1\ $
\begin{displaymath}
	\det(A_{ij} \cdot v_1)
	= (\alpha + h^{0,1}_{11}) (\alpha + h^{0,1}_{22})
	- h^{0,1}_{12} h^{0,1}_{21}
	= \alpha^2 + h^{0,1}_{11} h^{0,1}_{22}
	- h^{0,1}_{12} h^{0,1}_{21}
\end{displaymath}
when recalling that $\ g^{ij} A^0_{ij} = 0\ $.
Therefore
\begin{displaymath}
	| \det(A_{ij} \cdot v_1) - \alpha^2| \leq C |A^0|^2_g,
\end{displaymath}
which establishes the desired estimate with the equations above.
\proof
We are ready to prove Theorem \ref{intro.theorem-funda}.
\\ \ \\
{\pr Proof of Theorem \ref{intro.theorem-funda}:} \\
Theorem \ref{intro.theorem-funda} immediately follows
from Proposition \ref{glob-smooth.sphere}
for $\ \parallel A^0_\Sigma \parallel_{L^2(\Sigma)} \leq \delta_0(n) \ $
small enough depending on $\ n\ $.

For $\ \parallel A^0_\Sigma \parallel_{L^2(\Sigma)} \geq \delta_0(n) > 0\ $,
we conclude by the Gau\ss\ equations and the Gau\ss-Bonnet theorem
in (\ref{intro.gauss}) and (\ref{intro.gauss-impl}) that
\begin{displaymath}
	\parallel A_\Sigma \parallel_{L^2(\Sigma)}
	\leq \Cn \parallel A^0_\Sigma \parallel_{L^2(\Sigma)}.
\end{displaymath}
Then for any measureable unit normal vector field
$\ \nor \mbox{ on } \Sigma\ $
\begin{displaymath}
	\parallel A - \nor g \parallel_{L^2(\Sigma)}
	\leq \parallel A \parallel_{L^2(\Sigma)}
	+ 2 \Ht(\Sigma)^{1/2}
	\leq \Cn \parallel A^0_\Sigma \parallel_{L^2(\Sigma)}
\end{displaymath}
and, as $\ |K| \leq |A|^2/2\ $,
\begin{displaymath}
	\parallel K - 1 \parallel_{L^1(\Sigma)}
	\leq \frac{1}{2} \parallel A \parallel_{L^2(\Sigma)}^2
	+ \Ht(\Sigma)
	\leq \Cn \parallel A^0_\Sigma \parallel_{L^2(\Sigma)}^2,
\end{displaymath}
which finishes the proof also
for $\ \parallel A^0_\Sigma \parallel_{L^2(\Sigma)} \geq \delta_0(n)\ $.
\proof

%%%%%

%%%%%
%%%%%

\setcounter{equation}{0}

\section{Conformal parametrization} \label{conf}

In this section, we get a conformal parametrization
for surfaces in $ \ \rel^n \ $
extending \cite{del.muell} Proposition 3.2 to arbitrary codimension.
In order to formulate our result we recall the definition
of the number $\energ{n} = e_n$ from \cite{schae.comp-will12} (1.2).
\begin{displaymath}
	\energ{n}
	:= \left\{
	\begin{array}{cl}
		4 \pi & \mbox{for } n = 3, \\
		8 \pi / 3 & \mbox{for } n = 4, \\
		2 \pi &  \mbox{for } n \geq 5. \\
	\end{array}
	\right.
\end{displaymath}

\begin{proposition} \label{conf.bound}
Let $ \ \Sigma \subseteq {\mathbb R}^n \ $ be a smooth embedded surface of sphere type $ \ \Sigma \cong S^2 \ $ with 
\begin{equation}\label{conf.bound.area}
	\Ht(\Sigma) = 4 \pi 
\end{equation}
and
\begin{equation} \label{conf.bound.funda-0}
	\| A^0_\Sigma\|^2_{L^2(\Sigma)} < 2 \energ{n} .
\end{equation}
Then there exists a conformal parametrization $ \ f : S^2 \diff \Sigma \ $
with pull-back metric $\ g = f^* \geu = e^{2u}  g_{S^2}\ $ and
\begin{equation}\label{conf.bound.esti}
	\|u\|_{L^\infty(S^2)} \leq C(n,\tau),
\end{equation}
where $\tau =2 \energ{n} -\| A^0_\Sigma\|^2_{L^2(\Sigma)}$.
\end{proposition}
{\pr Proof:}\\
Since $\Sigma$ is a smooth spherical Riemann surface, the uniformization theorem (see \cite{far.kra}
Theorem IV. 4.1 or \cite{jo.rie}, Theorem 4.4.1) implies the existence of a conformal parametrization $f:S^2 \to \Sigma \subset \rel^n$.
Without loss of generality we assume that $f(e_3)=e_3$. Now we look at the inversion of $\rel^n$ at the sphere of radius $\sqrt{2}$ centered at $e_3$ which is given by 
\begin{displaymath}
\Phi(x) := e_3 + 2 (x - e_3) / |x - e_3|^2.
\end{displaymath}
Note that $\Phi | S^2: \partial B_1(0) \cap span \{ e_1 , e_2 , e_3 \}
\rightarrow \com \cup \{ \infty \}$ is the standard stereographic projection. 

We use $\Phi$ to define a conformal immersion $\hat f:\com \to \rel^n$ by 
\[
\hat f:= \Phi \circ f \circ \Phi^{-1}.
\]
The resulting image surface $\hat f(\com)=:\hat \Sigma \subset \rel^n$ is complete, connected and non-compact. Moreover, arguing as in the proof of Lemma 3.1 in \cite{kuw.schae.will3}, we get that
\begin{equation}
\int_{\hat \Sigma} \hat K_{\hat \Sigma} \d \Ht=0. \label{conf.bound.gauss.van}
\end{equation}
Note that here we deal with a smooth surface and therefore we don't need the assumption that the immersion is Willmore away from a possible singular point.  Next, using (\ref{conf.bound.funda-0}), (\ref{conf.bound.gauss.van}) and by pointwise conformal invariance of $\ |A^0_\Sigma|^2 \d \Ht\ $, see \bcite{chen.conf}, we estimate
\begin{equation}
\int_{\hat \Sigma} |\hat A_{\hat \Sigma}|^2  \d \Ht= 2 \int_{\hat \Sigma} |\hat A_{\hat \Sigma}^0|^2  \d \Ht
= 2 \int_\Sigma |A_{\Sigma}^0|^2  \d \Ht \le 4 \energ{n}- 2\tau . \label{conf.bound.est.secfun}
\end{equation}
The above facts allow us to apply Theorem 5.1 in \cite{schae.comp-will12} and hence we get the estimate
\begin{equation}
||\hat u-\hat \lambda||_{L^\infty(\com)} \le C(n,\tau) , \label{conf.bound.est.confhat}
\end{equation}
where $\hat \lambda =\lim_{z\to \infty} \hat u(z) \in \rel$ and $\hat f^* \geu =e^{2\hat u} g_\com$. Dilating in $\com$, which is a conformal parameter change for $f$, we may assume that $\hat \lambda = 0$.

Moreover, arguing as in the proof of Theorem 4.3.1 in \cite{muell.sver}
and combining with Theorem 5.1 in \cite{schae.comp-will12}, we get that
\begin{equation}
e^{-2C(n,\tau)}|z-w|\le |\hat f(z)-\hat f(w)|\le e^{C(n,\tau)}|z-w|. \label{conf.bound.bilipschitz}
\end{equation}
It follows from (\ref{conf.bound.funda-0}) that
\begin{displaymath}
4\pi \le \W(\Sigma) \le 8\pi-\tau/2.
\end{displaymath}
Combining this with (\ref{conf.bound.area}) and Lemma $1.1$ in \cite{sim.will} we obtain that the diameter of $\Sigma$ is bounded from below and above by
\begin{equation}
\sqrt{1/2} \le diam (\Sigma) \le C(n) \label{conf.bound.diam}
\end{equation}
where $C(n)$ is a constant which only depends on $n$. Therefore there exists another constant $C_1(n)>0$ so that for every $z\in \com$ we have
\begin{equation}
| \hat f(z) -e_3| \ge C_1(n)\label{conf.bound.est.e3}
\end{equation}
and, after a translation in $\com$, we can assume that
\begin{equation}
|\hat f(0) -e_3|\le 16. \label{conf.bound.est.e3a}
\end{equation}

We note that the derivative of $\Phi$ satisfies
\begin{displaymath}
\partial_i \Phi^j (x)=\frac{2}{|x-e_3|^2} \left(\delta_{ij}-\frac{2(x^i-e_3)(x^j-e_3)}{|x-e_3|^2} \right),
\end{displaymath}
where the matrix in the brackets is the orthogonal Householder matrix, and therefore the conformal factor $u_\Phi$ of $\Phi$ can be computed to be
\begin{displaymath}
u_\Phi(x) = \frac12 \log 4 -2 \log |x-e_3|.
\end{displaymath}
Since the conformal factor of a composition of conformal maps is the sum of the corresponding conformal factors, we get for every $x\in S^2\backslash \{e_3\}$ (note that $\Phi^{-1}=\Phi$)
\begin{equation}
u (x)  = \hat u(\Phi(x))  +2\log |\hat f(\Phi(x))-e_3|-2\log |\Phi(x)-e_3|,\label{conf.bound.eqn.confhat}
\end{equation}
where we also used that $f= \Phi^{-1} \circ \hat f \circ (\Phi| S^2)$.

Our first observation is that because of the estimate (\ref{conf.bound.bilipschitz}) the limit 
\begin{displaymath}
\lambda:=\lim_{x\to e_3} u(x)= \lim_{z\to \infty} (\hat u(z)+2\log |\hat f(z)-e_3|-2\log |z-e_3|)
\end{displaymath}
exists and we have
\begin{equation}
-2C(n,\tau)\le \lambda \le C(n,\tau).\label{conf.bound.est.uinfty}
\end{equation}

In the following we want to obtain uniform bounds from above and below for the quotient $|\hat f(z)-e_3| / |z-e_3|$, $z\in \com$.

In order to get these estimates, we note that (\ref{conf.bound.bilipschitz}) and (\ref{conf.bound.est.e3a}), imply that 
\begin{displaymath}
|\hat f(z)-e_3| \le |\hat f(z)-\hat f(0)|+|\hat f(0)-e_3| \le e^{C(n,\tau)}|z-e_3|+16+e^{C(n,\tau)}
\end{displaymath}
and
\begin{displaymath}
|\hat f(z)-e_3| \ge |\hat f(z)-\hat f(0)|-|\hat f(0)-e_3| \ge e^{-2C(n,\tau)}|z-e_3|-16-e^{-2C(n,\tau)}.
\end{displaymath}
Altogether, there exists a constant $C<\infty$ such that for every $z\in \com$ with $|z-e_3|\ge C(1+e^{2C(n,\tau)})$ we have
\begin{displaymath}
\frac12 e^{-2C(n,\tau)} |z-e_3| \le |\hat f(z)-e_3| \le 2 e^{C(n,\tau)} |z-e_3|.
\end{displaymath}
Hence we get that for every $z\in \com$ in the above range there exists a constant $C_1(n,\tau)>0$ so that
\begin{displaymath}
(C_1(n,\tau))^{-1} \le \log \left( \frac{|\hat f(z)-e_3|}{|z-e_3|} \right) \le C_1(n,\tau).
\end{displaymath}
For $1\le |z-e_3| =\sqrt{1+|z|^2} \le C(1+e^{2C(n,\tau)})$ we use (\ref{conf.bound.est.e3}) in order to get
\begin{displaymath}
C_1(n)(2C)^{-1}e^{-2C(n,\tau)} \le \frac{|\hat f(z)-e_3|}{|z-e_3|}\le 3e^{C(n,\tau)+\log 16}.
\end{displaymath}
Altogether this shows that for all $z\in \com$ there exists a constant $C_2(n,\tau)>0$ so that
\begin{displaymath}
(C_2(n,\tau))^{-1}\le \log \left( \frac{|\hat f(z)-e_3|}{|z-e_3|} \right) \le C_2(n,\tau).
\end{displaymath}
Inserting this estimate and (\ref{conf.bound.est.confhat}), (\ref{conf.bound.est.uinfty}) into the equation for $u$, we obtain that there exists a constant which we again call $C(n,\tau)$, so that
\begin{displaymath}
\|u\|_{L^\infty(S^2)} \le C(n,\tau).
\end{displaymath}
\proof
{\large \bf Remarks:} \\
\vspace{-.5cm}
\begin{enumerate}
\item
Assuming embeddedness in the above proposition is no restriction,
as we see by the Gau\ss\ equations
and the Gau\ss-Bonnet theorem in (\ref{intro.gauss})
and by (\ref{conf.bound.funda-0})
that $\ \W(\Sigma) < 8 \pi\ $,
since $\ \chi(\Sigma) \leq 2\ $,
and embeddedness follows from the Li-Yau inequality in \bcite{li.yau}
or \bcite{kuw.schae.will3} (A.17).

\item
(\ref{conf.bound.funda-0}) can equivalently be rewritten by (\ref{intro.gauss})
\begin{equation} \label{conf.bound.will}
	\W(\Sigma) < \energ{n} + 4 \pi =
	\left\{
	\begin{array}{cl}
		8 \pi & \mbox{for } n = 3, \\
		20 \pi / 3 & \mbox{for } n = 4, \\
		6 \pi &  \mbox{for } n \geq 5. \\
	\end{array}
	\right.
\end{equation}

\item
Actually for an embedded closed surface $\ \Sigma \subseteq \rel^n\ $
with either (\ref{conf.bound.funda-0}) and $\ \Sigma\ $ being orientable or the stronger condition that
\begin{equation} \label{conf.bound.funda-0-strong}
	\parallel A^0_\Sigma \parallel_{L^2(\Sigma)}^2 < 4 \pi,
\end{equation}
we see by (\ref{intro.gauss}), as $\ \W(f) \geq 4 \pi\ $
by the Li-Yau inequality in \bcite{li.yau}
or \bcite{kuw.schae.will3} (A.17),
that $\ \chi(\Sigma) > 0 \mbox{ respectively } \chi(\Sigma) > 1\ $,
hence $\ \chi(\Sigma) = 2\ $,
and we conclude $\ \Sigma \cong S^2\ $ is a sphere.
Clearly (\ref{conf.bound.funda-0-strong}) implies (\ref{conf.bound.funda-0}),
as $\ \energ{n} \geq 2 \pi\ $.

\end{enumerate}
\defin
\ \\
If $\ A^0\ $ is small in $\ L^2\ $ and $\ f\ $ is close to a round sphere in $\ L^\infty\ $,
we can prove smallness of the conformal factor in $\ L^\infty\ $.

\begin{proposition} \label{conf.small}

Let $\ f: S^2 \rightarrow \rel^n\ $ be a conformal immersion
with pull-back metric $\ g = f^*\geu = e^{2u} g_{S^2}\ $
for the canonical metric $\ g_{S^2} \mbox{ on } S^2\ $
satisfying
\begin{equation} \label{conf.small.area}
	area_g(S^2) = \int e^{2u} \d area_{S^2} = 4 \pi,
\end{equation}
\begin{equation} \label{conf.small.funda-0}
	\parallel A^0_f \parallel_{L^2(S^2,\mu_g)} \leq \delta
\end{equation}
and for $\ S^2 := \partial B_1(0)
\cap span \{ e_1 , e_2 , e_3 \} \subseteq \rel^n\ $
\begin{equation} \label{conf.small.round}
	\parallel f - id_{S^2} \parallel_{L^\infty(S^2)} \leq \delta.
\end{equation}
Then for $\ \delta \leq 1/2\ $
\begin{equation} \label{conf.small.conf-fak}
	\parallel u \parallel_{L^\infty(S^2)}
	\leq \Cn \delta.
\end{equation}
\end{proposition}
{\pr Proof:} \\
We continue in the notation of the previous Proposition \ref{conf.bound}
for the above $\ f\ $
and get as in (\ref{conf.bound.est.confhat})
by combining Theorem 5.1 in \cite{schae.comp-will12}
with (\ref{conf.bound.est.secfun}) and (\ref{conf.small.funda-0}) that
\begin{displaymath}
	osc_\com \hat u
	\leq \Cn \int \limits_{\hat \Sigma} |A_{\hat \Sigma}|^2 \d \Ht
	\leq \Cn \delta^2,
\end{displaymath}
as $\ \delta^2 \leq 1/4 \leq 2 \pi \leq e(n) / 2\ $.
We can rewrite the equation (\ref{conf.bound.eqn.confhat})
for every $x\in S^2\backslash \{e_3\}$ as follows
\begin{displaymath}
u(x) = \hat u(\Phi(x)) +2\log \frac{|f(x)-e_3|}{|x-e_3|}.
\end{displaymath}
Moreover for $\ x \in S_3^- := S^2 \cap \{ x_3 \leq 0 \}\ $, we have $\ |x - e_3| \geq 1\ $
and get from (\ref{conf.small.round}) that
\begin{displaymath}
	1-\delta \le \frac{|x-e_3| -|f(x)-x|}{|x-e_3|}
	\le \frac{|f(x)-e_3|}{|x-e_3|} \le \frac{|f(x)-x|+|x-e_3|}{|x-e_3|}
	\le 1+\delta,
\end{displaymath}
and hence conclude with $\ \delta \leq 1/2\ $ that 
\begin{displaymath}
	osc_{S^-_3} u \le \Cn \delta.
\end{displaymath}
After a rotation we can repeat the same argument in order to get
that the oscillation of $u$ on $S^+_3:= S^2 \cap \{ x_3 \geq 0 \}$
is also bounded by $C(\Lambda,n)\delta$. Altogether, this gives
\begin{displaymath}
osc_{S^2} u\le \Cn \delta.
\end{displaymath}
Finally we note that the assumption (\ref{conf.small.area}) implies the existence of a point $p\in S^2$ with $u(p)=0$ and therefore we conclude that
\begin{displaymath}
\parallel u\parallel_{L^\infty(S^2)} \le \Cn \delta,
\end{displaymath}
which is (\ref{conf.small.conf-fak}).
\proof
Next we turn to get a conformal parametrization
with equally distributed area on the half spheres
$\ S_i^\pm := S^2 \cap \{ \pm x_i \geq 0 \}, i = 1,2,3\ $.

\begin{proposition} \label{conf.half}

Let $\ g = e^{2u} g_{S^2}\ $ be a conformal metric on the sphere $\ S^2\ $ with
\begin{equation} \label{conf.half.area}
	area_g(S^2) = 4 \pi
\end{equation}
and
\begin{equation} \label{conf.half.conf-ass}
	\parallel u \parallel_{L^\infty(S^2)}
	\leq \Lambda
\end{equation}
for some $\ \Lambda < \infty\ $.
Then there exists a M\"obius transformation $\ \phi \mbox{ of } S^2\ $
such that the transformed metric
$\ \phi^* g = e^{2 v} g_{S^2}\ $ satisfies
\begin{equation} \label{conf.half.conf-bound}
	\parallel v \parallel_{L^\infty(S^2)}
	\leq C(\Lambda)
\end{equation}
and
\begin{equation} \label{conf.half.half}
	area_{\phi^* g}(S_i^\pm) = 2 \pi
\end{equation}
for the half spheres $\ S_i^\pm := S^2 \cap \{ \pm x_i \geq 0 \}, i = 1,2,3\ $.
\end{proposition}
{\pr Proof:} \\
We consider again the stereographic projection
$\ T: S^2 - \{ e_3 \} \diff \com\ $
and see
\begin{equation} \label{conf.half.radius}
	T: S_3^+ \diff \com - B_1(0),
	\quad
	T: S_3^- \diff \overline{B_1(0)}.
\end{equation}
We define
\begin{displaymath}
	area(r) := area_g(T^{-1}(B_r(0)))
\end{displaymath}
and see that $\ area\ $ is continuous,
since $\ area_g(T^{-1}(\partial B_r(0))) = 0\ $.
As clearly $\ area(r) \rightarrow 0 \mbox{ for } r \rightarrow 0
\mbox{ and } area(r) \rightarrow area_g(S^2) = 4 \pi
\mbox{ for } r \rightarrow \infty\ $,
there exists $\ 0 < r < \infty \mbox{ with } area(r) = 2 \pi\ $.
Recalling
\begin{displaymath}
	(T^{-1})^* g_{S^2} = \frac{4}{(1 + |.|^2)^2} \geu,
\end{displaymath}
we estimate with (\ref{conf.half.conf-ass}) that
\begin{displaymath}
	2 \pi = area_g(T^{-1}(B_r(0)))
	= \int \limits_{B_r(0)} e^{2 u (T^{-1}(z))}
	\frac{4}{(1 + |z|^2)^2} \d \Lt(z)
	\leq C(\Lambda) r^2,
\end{displaymath}
hence $\ r \geq c_0(\Lambda)\ $.
Inverting the complex plane, we obtain likewise
that $\ 1/r \geq c_0(\Lambda)\ $.
Dilating by $\ 1/r \mbox{ in } \com\ $
results in a M\"obius transformation $\ \phi \mbox{ on } S^2
\mbox{ with } \phi(x) := T^{-1}( T(x) / r )\ $
and by (\ref{conf.half.radius})
\begin{displaymath}
	area_{\phi^* g}(S_3^-)
	= area_{\phi^* g}(T^{-1}(\overline{B_1(0)}))
	= area_g(T^{-1}(\overline{B_r(0)})) = 2 \pi
\end{displaymath}
and $\ area_{\phi^* g}(S_3^+)
= area_{\phi^* g}(S^2) - area_{\phi^* g}(S_3^-) = 2 \pi\ $.
Moreover
\begin{displaymath}
	\phi^* g = e^{2 u \circ \phi} \phi^* g_{S^2}
	= e^{2 u \circ \phi} T^* (z \mapsto z/r)^* (T^{-1})^* g_{S^2}
	= e^{2 u \circ \phi} \frac{(1 + |T|^2)^2}{r^2 (1 + |T/r|^2)^2} g_{S^2},
\end{displaymath}
which yields (\ref{conf.half.conf-bound}) for $\ \phi^* g\ $
by (\ref{conf.half.conf-ass}) and by $\ c_0(\Lambda) \leq r \leq C(\Lambda)\ $.
Observing as in (\ref{conf.half.radius})
\begin{displaymath}
	T: S_i^\pm \diff \com \cap \{ \pm x_i \geq 0 \}
	\quad \mbox{for } i = 1,2,
\end{displaymath}
we see that $\ \phi: S_i^\pm \diff S_i^\pm \mbox{ for } i = 1,2,\ $
in particular $\ area_{\phi^* g}(S_i^\pm) = area_g(S_i^\pm)\ $.

Therefore applying the above procedure successively
to the stereographic projections at $\ e_3, e_2, \mbox{ and } e_1\ $,
we obtain (\ref{conf.half.half})
while keeping (\ref{conf.half.conf-bound})
with increasing constants.
This finishes the proof of the proposition.
\proof

%%%%%
%%%%%

\setcounter{equation}{0}
\section{Closeness to a round sphere} \label{round}

In this section, we consider a conformal immersion
$\ f: S^2 \rightarrow \rel^n\ $
with pull-back metric $\ g = f^*\geu = e^{2u} g_{S^2}\ $
for the canonical metric $\ g_{S^2} \mbox{ on } S^2\ $
and induced measure $\ \mu_g = \sqrt{g}\ area_{S^2}\ $ satisfying
\begin{equation} \label{round.ass.conf-bound}
	\parallel u \parallel_{L^\infty(S^2)}
%	\parallel u \parallel_{W^{1,2}(S^2)}
	\leq \Lambda
\end{equation}
for some $\ \Lambda < \infty\ $,
\begin{equation} \label{round.ass.funda-0}
	\parallel A^0_f \parallel_{L^2(S^2,\mu_g)} \leq \delta
\end{equation}
normalized by
\begin{equation} \label{round.ass.area}
	area_g(S^2) = \int e^{2u} \d area_{S^2} = 4 \pi
\end{equation}
and nearly equally distributed area on the half spheres 
$\ S_i^\pm := S^2 \cap \{ \pm x_i \geq 0 \}, i = 1,2,3,\ $
\begin{equation} \label{round.ass.half}
	| area_g(S_i^\pm) - 2 \pi | \leq \delta
\end{equation}
for some $\ \delta > 0\ $ small.

First we estimate the exponent of the conformal factor in $\ L^2\ $.

\begin{proposition} \label{round.conf}

Under the assumptions (\ref{round.ass.conf-bound})
- (\ref{round.ass.half}), we get
\begin{equation} \label{round.conf.conf-small}
	\parallel u \parallel_{L^2(S^2)} \leq C(\Lambda,n) \delta
\end{equation}
for $\ \delta\ $ small enough depending on $\ \Lambda,n\ $.
\end{proposition}
{\pr Proof:} \\
By the Gauss-equations, we know
\begin{equation} \label{round.conf.gauss-equ}
	-\Delta_{S^2} u = K e^{2u} - 1.
\end{equation}
By Proposition \ref{glob-smooth.sphere} (\ref{glob-smooth.sphere.gauss}),
we have
\begin{displaymath}
	\parallel K - 1 \parallel_{L^1(S^2)}
	\leq C(\Lambda,n) \delta
\end{displaymath}
and see that $\ K \rightarrow 1 \mbox{ strongly in } L^1(S^2)
\mbox{ when } \delta \rightarrow 0\ $.
Mulitplying (\ref{round.conf.gauss-equ}) by $\ u\ $, we get
\begin{displaymath}
	\int \limits_{S^2} |\nabla u|^2 \d area_{S^2}
	= \int \limits_{S^2} (K e^{2u} - 1) u \d area_{S^2}
	\leq C(\Lambda),
\end{displaymath}
hence together with (\ref{round.ass.conf-bound})
\begin{displaymath}
	\parallel u \parallel_{W^{1,2}(S^2)} \leq C(\Lambda).
\end{displaymath}
We follow \bcite{del.muell} \S 6.1 and first get
\begin{equation} \label{round.conf.conf-small-aux}
	\parallel u \parallel_{L^2(S^2)} \leq \eta
\end{equation}
for any $\ \eta > 0\ $, if $\ \delta\ $ is small enough.
Writing (\ref{round.conf.gauss-equ})
as in \bcite{del.muell} \S 6.1 (75) in the form
\begin{displaymath}
	-\Delta_{S^2} u - 2 u
	= e^{2u} - 2u - 1 + (K-1) e^{2u},
\end{displaymath}
we get
\begin{equation} \label{round.conf.aux}
	\parallel \Delta_{S^2} u + 2 u \parallel_{L^1(S^2)}
	\leq C(\Lambda) \parallel u \parallel_{L^2(S^2)}^2
	+ C(\Lambda) \delta.
\end{equation}
Combining (\ref{round.ass.half}) and (\ref{round.conf.aux}),
we get from \bcite{del.muell} \S 6.1 (68)
\begin{displaymath}
	\parallel u \parallel_{L^2(S^2)}
	\leq C(\Lambda) \parallel u \parallel_{L^2(S^2)}^2
	+ C(\Lambda) \delta,
\end{displaymath}
which yields (\ref{round.conf.conf-small}),
when $\ C \parallel u \parallel_{L^2(S^2)} \leq C \eta \leq 1/2\ $,
which is true for $\ \delta\ $ small enough.
\proof
Now we a ready to prove that the immersion $\ f\ $
is close to a round sphere.

\begin{proposition} \label{round.sphere}

Under the assumptions (\ref{round.ass.conf-bound})
- (\ref{round.ass.half}),
we get after an appropriate translation and rotation
and with $\ S^2 := \partial B_1(0)
\cap span \{ e_1 , e_2 , e_3 \} \subseteq \rel^n\ $
\begin{equation} \label{round.sphere.esti}
	\parallel f - id_{S^2} \parallel_{W^{2,2}(S^2)}
	+ \parallel u \parallel_{L^\infty(S^2)}
	\leq C(\Lambda,n) \delta
\end{equation}
for $\ \delta\ $ small enough depending on $\ \Lambda,n\ $.
\end{proposition}
{\pr Proof:} \\
We see from $\ \Delta_g f = \mean_f\ $ that
\begin{displaymath}
	\Delta_{S^2} f + 2 f
	= e^{2u} (\mean_f + 2f) + (1 - e^{2u}) 2 f
\end{displaymath}
and estimate using
(\ref{glob-smooth.sphere.mean}),
(\ref{glob-smooth.sphere.l-infty}),
(\ref{round.ass.conf-bound}),
and (\ref{round.conf.conf-small}) that
\begin{displaymath}
	\parallel \Delta_{S^2} f + 2 f \parallel_{L^2(S^2)}
	\leq C(\Lambda,n) \delta.
\end{displaymath}
when recalling that $\ |e^{2u} - 1| \leq C(\Lambda) |u|\ $.
Since the kernel of $\ \Delta_{S^2} + 2\ $
consists of exactly the linear functions,
there exists a linear function $\ l = (l_1,l_2,l_3):
\rel^3 \supseteq S^2 \rightarrow \rel^n\ $ with
\begin{displaymath}
	\parallel f - l \parallel_{L^2(S^2)}
	\leq C \parallel \Delta_{S^2} f + 2 f \parallel_{L^2(S^2)}
	\leq C(\Lambda,n) \delta.
\end{displaymath}
see also \bcite{del.muell} \S 6.1 (70).
%%%%%
Likewise let $\ v_m\ $ be an orthonormal basis of $\ L^2(S^2)\ $
of eigenfunction of $\ -\Delta_{S^2} \mbox{ with eigenvalue }
\lambda_m \mbox{ for } m \in \nat_0\ $
and with $\ 0 = \lambda_0 < \lambda_1 \leq ... \leq \lambda_m \leq \ldots\ $
In particular $\ v_0\ $ is a constant
and $\ v_1, v_2, v_3\ $ form an orthonormal basis
of the linear functions in $\ \rel^3\ $
with $\ 2 = \lambda_1 = \lambda_2 = \lambda_3 < \lambda_4 \leq \ldots\ $
Putting $\ \alpha_m = \langle f , v_m \rangle_{S^2}\ $,
we see $\ f = \sum_{m=0}^\infty \alpha_m v_m \mbox{ in } L^2(S^2)\ $.
Moreover $\ 0 = \langle f , \Delta_{S^2} v_m + \lambda_m v_m \rangle_{S^2}
= \langle \Delta_{S^2} f + \lambda_m f , v_m \rangle_{S^2}\ $,
hence $\ (2 - \lambda_m) \alpha_m
= \langle \Delta_{S^2} f + 2 f , v_m \rangle_{S^2}\ $.
Now $\ l = \alpha_1 v_1 + \alpha_2 v_2 + \alpha_3 v_3\ $ is linear,
and we estimate
\begin{displaymath}
	\parallel f - l \parallel_{L^2(S^2)}^2
	= \sum \limits_{m=0, m \neq 1,2,3}^\infty |\alpha_m|^2
	= \sum \limits_{m=0, m \neq 1,2,3}^\infty |2 - \lambda_m|^{-2}
	\ |\langle \Delta_{S^2} f + 2 f , v_m \rangle_{S^2}|^2 \leq
\end{displaymath}
\begin{displaymath}
	\leq \min(2, |\lambda_4 - 2|)^{-2}
	\parallel \Delta_{S^2} f + 2 f \parallel_{L^2(S^2)}^2,
\end{displaymath}
as desired.

%%%%%
%\input{ch-round-sphere-l1}
As $\ \Delta_{S^2} l + 2 l = 0\ $,
we get by standard elliptic theory,
see \bcite{gil.tru} Theorem 8.8, that
\begin{equation} \label{round.sphere.esti-aux}
	\parallel f - l \parallel_{W^{2,2}(S^2)}
	\leq \Cn \parallel \Delta_{S^2} f + 2 f \parallel_{L^2(S^2)}
	+ \Cn \parallel f - l \parallel_{L^2(S^2)}
	\leq C(\Lambda,n) \delta.
\end{equation}
In particular there exists $\ p \in S^2\ $ with
\begin{equation} \label{round.sphere.conf-aux}
	| \nabla f(p) - l | \leq C(\Lambda,n) \delta
	\quad \mbox{and} \quad
	|e^{2u(p)} - 1| \leq C(\Lambda,n) \delta.
\end{equation}
Assuming after a rotation that $\ p = e_3\ $,
we have $\ T_p S^2 = span \{ e_1 , e_2 \}\ $ and conclude,
as $\ \partial_i f \cdot \partial_j f = e^{2u} \delta_{ij}\ $, that
\begin{displaymath}
	|\langle l_i , l_j \rangle - \delta_{ij} |
	\leq |\langle l_i , l_j \rangle - e^{2u(p)} \delta_{ij} |
	+ |e^{2u(p)} - 1|
	\leq C(\Lambda,n) \delta
	\quad \mbox{for } i,j = 1,2.
\end{displaymath}
Applying a Gram-Schmidt orthonormalisation to $\ l_1, l_2\ $,
we obtain orthonormal vectors $\ \tilde l_1, \tilde l_2
\mbox{ with } |\tilde l_1 - l_1|, |\tilde l_2 - l_2|
\leq C(\Lambda,n) \delta \mbox{ for } \delta\ $ small enough.
Replacing $\ l_1, l_2 \mbox{ by } \tilde l_1, \tilde l_2\ $,
we may assume that $\ l_1, l_2\ $ are orthonormal
while keeping the estimate (\ref{round.sphere.esti-aux}).

Next there exists $\ q \in S^2 \mbox{ with } |q_3| \leq 1/\sqrt{2}\ $
satisfying (\ref{round.sphere.conf-aux})
with $\ p \mbox{ replaced by } q\ $.
After a rotation, we may assume that $\ q_1 = 0\ $.
Then $\ T_q S^2 = span \{ e_1 , (0,-q_3,q_2) \}\ $,
in particular $\ T_q S^2 \cap span \{ e_1 , e_2 \} = span \{ e_1 \}\ $,
and as above
\begin{displaymath}
	| \langle l_1 , -q_3 l_2 + q_2 l_3 \rangle |
	\leq C(\Lambda,n) \delta,
\end{displaymath}
hence, as $\ l_1 \perp l_2\ $, that
\begin{equation} \label{round.sphere.ortho-aux}
	| \langle l_1 , l_3 \rangle |
	\leq C(\Lambda,n) \delta / |q_2|	
	\leq C(\Lambda,n) \delta,
\end{equation}
since $\ |q_2| = \sqrt{1 - |q_3|^2} \geq 1/\sqrt{2}\ $.

In the same way there exists $\ r \in S^2 \mbox{ with } |r_3| \leq 1/\sqrt{2}
\mbox{ and } r_1 \geq 1/2\ $
satisfying (\ref{round.sphere.conf-aux})
with $\ p \mbox{ replaced by } r\ $.
Then $\ span \{ (-r_2,r_1,0) \} = T_r S^2 \cap span \{ e_1 , e_2 \} \neq \{ 0 \}\ $
and with the same argument as above, we obtain
\begin{displaymath}
	| \langle -r_2 l_1 + r_1 l_2 , l_3 \rangle |
	\leq C(\Lambda,n) \delta\ |(-r_2,r_1,0)|
	\leq C(\Lambda,n) \delta,
\end{displaymath}
hence combining with (\ref{round.sphere.ortho-aux})
\begin{displaymath}
	|\langle l_2 , l_3 \rangle|
	\leq r_1^{-1} \Big( C(\Lambda,n) \delta
	+ |r_2|\ |\langle l_1 , l_3 \rangle| \Big)
	\leq C(\Lambda,n) \delta,
\end{displaymath}
since $\ r_1 \geq 1/2, |r_2| \leq 1\ $.

Returning to $\ q \in S^2\ $, we further obtain
\begin{displaymath}
	C(\Lambda,n) \delta \geq |\ |-q_3 l_2 + q_2 l_3|^2 - 1|
	= |1 - q_3^2 - q_2^2 |l_3|^2 + 2 q_2 q_3 \langle l_2 , l_3 \rangle\ |
	\geq q_2^2\ |1 - |l_3|^2| - C(\Lambda,n) \delta 
\end{displaymath}
hence
\begin{displaymath}
	|\ |l_3| - 1| \leq C(\Lambda,n) \delta.
\end{displaymath}
Now applying a Gram-Schmidt orthonormalisation to $\ l_1, l_2, l_3\ $,
we obtain orthonormal vectors $\ l_1, l_2, \tilde l_3
\mbox{ with } |\tilde l_3 - l_3| \leq C(\Lambda,n) \delta
\mbox{ for } \delta\ $ small enough,
and replacing $\ l_3 \mbox{ by } \tilde l_3\ $,
we may assume that $\ l_1, l_2, l_3\ $ are orthonormal
while keeping the estimate (\ref{round.sphere.esti-aux}).

After a further rotation, we may assume that $\ l = (e_1,e_2,e_3)\ $,
hence $\ l = id_{S^2} \mbox{ on } S^2 
= \partial B_1(0) \cap span \{ e_1 , e_2 , e_3 \} \subseteq \rel^n\ $,
and the first estimate in (\ref{round.sphere.esti})
follows directly from (\ref{round.sphere.esti-aux}).

For the $\ L^\infty$-estimate for $\ u\ $ in (\ref{round.sphere.esti}),
we note that the Sobolev-embedding
$\ W^{2,2}(B_1^2(0)) \hookrightarrow L^\infty(B_1^2(0))\ $ implies
\begin{displaymath}
	\parallel f - id_{S^2} \parallel_{L^\infty(S^2)}
	\leq C(\Lambda,n) \delta
\end{displaymath}
and obtain by Proposition \ref{conf.small}
for $\ C(\Lambda,n) \delta \leq 1/2\ $ that
\begin{displaymath}
	\parallel u \parallel_{L^\infty(S^2)}
	\leq C(\Lambda,n) \delta,
\end{displaymath}
which establishes the second estimate in (\ref{round.sphere.esti})
and concludes the proof of the proposition.
\proof
Finally we prove Theorem \ref{intro.theorem-round}.
\\ \ \\
{\pr Proof of Theorem \ref{intro.theorem-round}:} \\
First we get by Propositions \ref{conf.bound} and \ref{conf.half}
a conformal parametrization $\ f: S^2 \diff \Sigma\ $ satisfying
(\ref{round.ass.conf-bound}) - (\ref{round.ass.half})
with $\ \Lambda = C(n,\tau)\ $ in (\ref{round.ass.conf-bound}),
$\ \delta = \parallel A^0_\Sigma \parallel_{L^2(\Sigma)}\ $
in (\ref{round.ass.funda-0})
and equality in (\ref{round.ass.half}).

For $\ \parallel A^0_\Sigma \parallel_{L^2(\Sigma)} \leq \delta_0(n) \ $
small enough depending on $\ n\ $,
we see $\ \tau = 2 \energ{n} - \delta^2 \geq 2 \pi\ $ say,
hence $\ \Lambda = \Cn\ $ and get by Proposition \ref{round.sphere}
after an appropriate translation and rotation
\begin{displaymath}
	\parallel f - id_{S^2} \parallel_{W^{2,2}(S^2)}
	+ \parallel u \parallel_{L^\infty(S^2)}
	\leq \Cn \parallel A^0_\Sigma \parallel_{L^2(\Sigma)},
\end{displaymath}
which establishes the theorem for
$\ \parallel A^0_\Sigma \parallel_{L^2(\Sigma)} \leq \delta_0(n)\ $.

For $\ \parallel A^0_\Sigma \parallel_{L^2(\Sigma)} \geq \delta_0(n) > 0\ $,
we see by the Gau\ss\ equations and the Gau\ss-Bonnet theorem
in (\ref{intro.gauss}) and (\ref{intro.gauss-impl}) that
\begin{displaymath}
%	\W(\Sigma),
	\parallel A_\Sigma \parallel_{L^2(\Sigma)}
	\leq \Cn \parallel A^0_\Sigma \parallel_{L^2(\Sigma)}.
\end{displaymath}
Next by $\ \Delta_{S^2} f = e^{2u} \mean_f
\mbox{ for } f^* \geu =: e^{2u} g_{S^2}\ $,
we get by standard elliptic theory,
see \bcite{gil.tru} Theorem 8.8, that
\begin{displaymath}
	\parallel f \parallel_{W^{2,2}(S^2)}
	\leq C(n,\tau) \parallel \mean_f \parallel_{L^2(S^2)}
	\leq C(n,\tau) \parallel A^0_\Sigma \parallel_{L^2(\Sigma)},
\end{displaymath}
hence
\begin{displaymath}
	\parallel f - id_{S^2} \parallel_{W^{2,2}(S^2)}
	\leq \parallel f \parallel_{W^{2,2}(S^2)}
	+ \parallel id_{S^2} \parallel_{W^{2,2}(S^2)} \leq
\end{displaymath}
\begin{displaymath}
	\leq C(n,\tau) \parallel A^0_\Sigma \parallel_{L^2(\Sigma)} + C
	\leq C(n,\tau) \parallel A^0_\Sigma \parallel_{L^2(\Sigma)}.
\end{displaymath}
Clearly (\ref{round.ass.conf-bound}) with $\ \Lambda = C(n,\tau)\ $ implies
\begin{displaymath}
	\parallel u \parallel_{L^\infty(S^2)}
	\leq C(n,\tau) \delta_0(n)^{-1}
	\parallel A^0_\Sigma \parallel_{L^2(\Sigma)},
\end{displaymath}
which finishes the proof also
for $\ \parallel A^0_\Sigma \parallel_{L^2(\Sigma)} \geq \delta_0(n)\ $.
\proof

%%%%%

%\newpage

%%%%%
\ \\
{\LARGE \bf Appendix}
%\\ \ \\
%In this appendix, we collect for the reader's convenience
%some results which are consequences or adaptions of standard
%results.

\begin{appendix}
\renewcommand{\theequation}{\mbox{\Alph{section}.\arabic{equation}}}

%%%%%

\setcounter{equation}{0}
\section{Hausdorff convergence} \label{haus}

\begin{proposition} \label{haus.meas-low}

Let $\ \mu\ $ be an integral $\ 2-$varifold
in an open set $\ U \subseteq \rel^n\ $
with square integrable weak mean curvature
$\ \mean_\mu \in L^2(\mu),
\parallel \mean_\mu \parallel_{L^2(\mu)}^2 \leq \W < \infty\ $
and connected support $\ spt\ \mu\ $.
Then for any $\ x_0 \in spt\ \mu
\mbox{ with } spt\ \mu \not\subseteq B_\varrho(x_0) \subseteq U\ $,
we have
\begin{equation} \label{haus.meas-low.esti}
	\mu(B_\varrho(x_0))
	\geq c_0 \varrho^2 / (1 + \W)
\end{equation}
for some $\ c_0 > 0\ $.
\end{proposition}
{\pr Proof:} \\
We select an integer $\ N\ $
and see, as $\ x_0 \in spt\ \mu \not\subseteq B_\varrho(x_0)
\mbox{ and } spt\ \mu\ $ is connected, that there exist
\begin{displaymath}
	x_m \in spt\ \mu \cap \partial B_{(j + 1/2) \varrho / N}(x_0)
	\neq \emptyset
	\quad \mbox{for } j = 1, \ldots, N-1.
\end{displaymath}
By \bcite{kuw.schae.will3} (A.6), (A.10),
\begin{displaymath}
	1 \leq \theta^2(\mu,x_m)
	\leq C (\varrho/(2N))^{-2} \mu(B_{\varrho/(2N)}(x_m))
	+ C \parallel \mean_\mu \parallel_{L^2(B_{\varrho/(2N)}(x_m))}^2
\end{displaymath}
and adding up observing that
$\ B_{\varrho/(2N)}(x_m), j = 0, \ldots, N-1\ $ are pairwise disjoint,
\begin{displaymath}
	N \leq C N^2 \varrho^{-2} \mu(B_\varrho(x_0)) + C \W.
\end{displaymath}
If $\ \W / (\varrho^{-2} \mu(B_\varrho(x_0))) \leq 1\ $,
we have (\ref{haus.meas-low.esti}) with $\ c_0 = 1 / (2 C)\ $.
Otherwise we can select an integer
$\ N \leq \sqrt{\W / (\varrho^{-2} \mu(B_\varrho(x_0)))}
\leq N + 1 \leq 2N\ $. This yields
\begin{displaymath}
	1 \leq C \sqrt{\W \varrho^{-2} \mu(B_\varrho(x_0))}
\end{displaymath}
which implies (\ref{haus.meas-low.esti}).
\proof

\begin{proposition} \label{haus.prop}

Let $\ \mu_m\ $ be a sequence of integral $\ 2-$varifolds
in an open set $\ U \subseteq \rel^n
\mbox{ with } \parallel \mean_{\mu_m} \parallel_{L^2(\mu_m)}^2
\leq \W < \infty $,
$\ \mu_m \rightarrow \mu\ $ weakly as varifolds.

If
\begin{equation} \label{haus.prop.ass-1}
	spt\ \mu_m \mbox{ are connected and }
	diam(spt\ \mu_m) \not\rightarrow 0
\end{equation}
or
\begin{equation} \label{haus.prop.ass-2}
	\W \leq \varepsilon_0
\end{equation}
for some $\ \varepsilon_0 > 0\ $ small enough,
then $\ spt\ \mu_m \rightarrow spt\ \mu\ $
locally in Hausdorff distance,
that is
\begin{equation} \label{haus.prop.haus}
	spt\ \mu = \{ x \in U\ |
	\ \exists x_m \in spt\ \mu_m:
	x_m \rightarrow x\ \}.
\end{equation}
If $\ spt\ \mu_m\ $ are connected
and for some open sets $\ U_1, U_2 \subseteq U\ $
\begin{equation} \label{haus.prop.ass-conn}
	spt\ \mu \subseteq U_1 \cup U_2,
	U_1 \cap U_2 = \emptyset,
	U_1 \subset \subset U,
	U_1 \cap spt\ \mu \neq \emptyset,
\end{equation}
then
\begin{equation} \label{haus.prop.aux}
	spt\ \mu \subseteq U_1
\end{equation}
and in particular
\begin{equation} \label{haus.prop.conn}
	spt\ \mu \mbox{ is connected}.
\end{equation}
\end{proposition}
{\pr Proof:} \\
For $\ x_0 \in spt\ \mu\ $,
we know for any $\ \varrho > 0\ $ that
$\ 0 < \mu(B_\varrho(x_0))
\leq \liminf_{m \rightarrow \infty}
\mu_m(B_\varrho(x_0))\ $,
hence $\ B_\varrho(x_0) \cap spt\ \mu_m
\neq \emptyset \mbox{ for large } m\ $.
By appropriate choice,
we find $\ x_m \in spt\ \mu_m
\mbox{ with } x_m \rightarrow x_0\ $.

Next consider $\ x_m \in spt\ \mu_m
\mbox{ with } x_m \rightarrow x_0 \in U\ $.
Under assumption (\ref{haus.prop.ass-1}),
we see $\ spt\ \mu_m \not\subseteq B_\varrho(x_m)\ $
for some $\ \varrho > 0 \mbox{ and large } m\ $.
As $\ spt\ \mu_m\ $ is further connected,
we get by (\ref{haus.meas-low.esti})
for any $\ 0 < r < \varrho\ $ that
\begin{displaymath}
	\mu_m(B_r(x_m)) \geq c_0 r^2 / (1 + \W),
\end{displaymath}
hence by weak convergence $\ \mu_m \rightarrow \mu\ $ that
\begin{displaymath}
	\mu(\overline{B_r(x_0)})
	\geq \limsup \limits_{m \rightarrow \infty}
	\mu_m(B_r(x_m)) \geq c_0 r^2 / (1 + \W)
\end{displaymath}
and $\ x_0 \in spt\ \mu\ $.

Assuming (\ref{haus.prop.ass-2}),
we consider again $\ x_m \in spt\ \mu_m
\mbox{ with } x_m \rightarrow x_0 \in U\ $.
By \bcite{kuw.schae.will3} (A.6) with $\ \delta = 1\ $ and (A.10),
we get for $\ x_m \in B_\varrho(x_0)
\subseteq B_{2 \varrho}(x_0) \subseteq U \mbox{ and } 0 < r < \varrho\ $
\begin{displaymath}
	r^{-2} \mu_m(B_r(x_m))
	\geq \omega_2/2  - C \W
	\geq \omega_2/3
\end{displaymath}
when $\ \W \leq \varepsilon_0\ $ is small enough,
hence again
\begin{displaymath}
	\mu(\overline{B_r(x_0)})
	\geq \limsup \limits_{m \rightarrow \infty}
	\mu_m(B_r(x_m)) \geq \omega_2 r^2 / 3
\end{displaymath}
and $\ x_0 \in spt\ \mu\ $.

Finally, we assume (\ref{haus.prop.ass-conn})
and that $\ spt\ \mu_m\ $ are connected.
If $\ spt\ \mu \cap U_2 \neq \emptyset\ $,
then by (\ref{haus.prop.haus}), we see
$\ spt\ \mu_m \cap U_i \neq \emptyset
\mbox{ for } m\ $ large enough.
Since $\ spt\ \mu_m\ $ is connected,
there exists $\ x_m \in spt\ \mu_m \cap \partial U_1
\neq \emptyset\ $.
Since $\ U_1 \subset \subset U\ $,
we get for a subsequence $\ x_m \rightarrow x_0
\in \partial U_1 \subseteq U\ $.
Again by (\ref{haus.prop.haus}),
we see $\ x_0 \in spt\ \mu\ $.
As $\ spt\ \mu \cap \partial U_1 = \emptyset\ $,
this is a contradiction,
hence $\ spt\ \mu \cap U_2 = \emptyset
\mbox{ and } spt\ \mu \subseteq U_1\ $,
which is (\ref{haus.prop.aux}).

If $\ spt\ \mu\ $ is not connected,
there are $\ V_1, V_2 \subseteq U_1\ $ open
with $\ spt\ \mu \subseteq V_1 \cup V_2,
V_1 \cap V_2 = \emptyset
\mbox{ and } spt\ \mu \cap V_i \neq \emptyset\ $.
Since $\ V_1 \subseteq U_1 \subset \subset U\ $,
we get by (\ref{haus.prop.aux})
that $\ spt\ \mu \cap V_2 = \emptyset\ $,
which is a contradiction,
and hence $\ spt\ \mu\ $ is connected.
\proof
{\large \bf Remark:} \\
Here we collect some consequences of the monotonicity
formula following \bcite{kuw.schae.will3} \S A.
Let $\ \mu\ $ be an integral $\ 2-$varifold in $\ \rel^n\ $
with square integrable weak mean curvature and
$\ \W(\mu) := (1/4) \parallel \mean_\mu \parallel_{L^2(\mu)}^2\ $.
%\vspace{-.5cm}
\begin{enumerate}
\item
\bcite{kuw.schae.will3} (A.14) says that
\begin{displaymath}
	\theta^2(\mu,\infty)
	:= \lim \limits_{\varrho \rightarrow 0}
	\frac {\mu(B_{\varrho}(0))}
	{\omega_2 \varrho^2} \in [0,\infty]
\end{displaymath}
exists,
in particular
\begin{displaymath}
	\frac{(\zeta_{R,\#} \mu)(B_\varrho(0))}{\omega_2 \varrho^2}
	= \frac{\mu(B_{R \varrho}(0))}{\omega_2 (R \varrho)^2}
	\rightarrow \theta^2(\mu,\infty)
	\quad \mbox{for all } \varrho > 0,
	\mbox{ as } R \rightarrow \infty.
\end{displaymath}
Passing to the limit $\ \zeta_{R_m,\#} \mu
\rightarrow \nu\ $ weakly as varifolds
for subsequences $\ R_m \rightarrow \infty\ $,
we see that $\ \nu\ $ is a stationary integral varifold,
as $\ \W(\mu) < \infty\ $, with
\begin{displaymath}
	\nu(B_\varrho(0)) = \theta^2(\mu,\infty) \omega_2 \varrho^2
	\quad \mbox{for all } \varrho > 0.
\end{displaymath}
By the argument as in \bcite{sim} \S 19 and \S 42,
we see that $\ \nu\ $ is a cone at the origin.
Inverting by $\ I(x) := x / |x|^2\ $,
we see that $\ I_\# \nu = \nu\ $.
Putting $\ \hat \mu := I_\# \mu\ $,
we see that
\begin{displaymath}
	\zeta_{\varrho_m^{-1},\#} \hat \mu
	= I_\# \zeta_{\varrho_m,\#} \mu
	\rightarrow I_\# \nu = \nu
	\quad \mbox{weakly as varifolds},
\end{displaymath}
hence
%\begin{displaymath}
%	\lim \limits_{m \rightarrow \infty}
%	\frac{\hat \mu(B_{\varrho_m^{-1} \varrho}(0))}
%	{\omega_2 (\varrho_m^{-1} \varrho)^2}
%	\leq \lim \limits_{m \rightarrow \infty}
%	\frac{(\zeta_{\varrho_m^{-1},\#} \hat \mu)(B_\varrho(0))}
%	{\omega_2 \varrho^2}
%	\leq \frac{\nu(B_\varrho(0))}{\omega_2 \varrho^2} \leq
%\end{displaymath}
%\begin{displaymath}
%	\leq \frac{\nu(\overline{B_\varrho(0)})}{\omega_2 \varrho^2}
%	\leq \lim \limits_{m \rightarrow \infty}
%	\frac{(\zeta_{\varrho_m^{-1},\#} \hat \mu)
%	(\overline{B_\varrho(0)})}
%	{\omega_2 \varrho^2}
%	\leq \lim \limits_{m \rightarrow \infty}
%	\frac{\hat \mu(\overline{B_{\varrho_m^{-1} \varrho}(0)})}
%	{\omega_2 (\varrho_m^{-1} \varrho)^2},
%\end{displaymath}
\begin{displaymath}
	\lim \limits_{m \rightarrow \infty}
	\frac{\hat \mu(B_{\varrho_m^{-1} \varrho}(0))}
	{\omega_2 (\varrho_m^{-1} \varrho)^2}
	= \frac{\nu(B_\varrho(0))}{\omega_2 \varrho^2}
	= \theta^2(\mu,\infty)
\end{displaymath}
and
\begin{equation} \label{haus.mono.inv}
	\theta^2(I_\# \mu,0) = \theta^2(\mu,\infty).
\end{equation}

\item
Combining \bcite{kuw.schae.will3} (A.8), (A.14)
and the inequality before (A.23),
we get for $\ \mu_m \mbox{ with } \W(\mu_m) < \infty
\mbox{ and } \mu_m \rightarrow \mu\ $ weakly as varifolds
\begin{equation} \label{haus.mono.infty}
	\theta^2(\mu,\cdot)
	\leq \theta^2(\mu,\infty)
	+ \frac{1}{4 \pi} \W(\mu)
	\leq \liminf \limits_{m \rightarrow \infty}
	\Big( \theta^2(\mu_m,\infty)
 	+ \frac{1}{4 \pi} \W(\mu_m) \Big).
\end{equation}
Actually putting $\ \beta_m := (1/4) |\mean_{\mu_m}|^2 \mu_m\ $
and passing to the limit for a subsequence
$\ \beta_m \rightarrow \beta\ $ weakly as Radon measures,
we can improve to
\begin{displaymath}
	\theta^2(\mu,\cdot)
	\leq \theta^2(\mu,\infty)
	+ \frac{1}{4 \pi} \W(\mu)
	+ \frac{1}{4} \sum \limits_{x \in \rel^n} \beta(\{x\}) \leq
\end{displaymath}
\begin{equation} \label{haus.mono.infty-sharp}
	\leq \theta^2(\mu,\infty)
	+ \frac{1}{4} \beta(\rel^n)
	\leq \liminf \limits_{m \rightarrow \infty}
	\Big( \theta^2(\mu_m,\infty)
 	+ \frac{1}{4 \pi} \W(\mu_m) \Big).
\end{equation}

\item
If $\ spt\ \mu\ $ consists
of at least $\ N\ $ components,
we can find in successive stages
open pairwise disjoint subsets
$\ U_1, \ldots, U_N \subseteq \rel^n\ $ with
\begin{displaymath}
\begin{array}{c}
	spt\ \mu \subseteq U_1 \cup \ldots \cup U_N, \\
	U_i \cap spt\ \mu \neq \emptyset
	\quad \mbox{for all } i = 1, \ldots, N,
\end{array}
\end{displaymath}
hence by (\ref{haus.mono.infty}) and \bcite{kuw.schae.will3} (A.10)
\begin{displaymath}
	N \leq \sum \limits_{i=1}^N
	\Big( \theta^2(\mu \lfloor U_i,\infty)
	+ \frac{1}{4 \pi} \W(\mu \lfloor U_i) \Big)
	= \theta^2(\mu,\infty)
	+ \frac{1}{4 \pi} \W(\mu)
\end{displaymath}
and
\begin{equation} \label{haus.mono.conn}
	\#( \mbox{components of } spt\ \mu )
	\leq \theta^2(\mu,\infty)
	+ \frac{1}{4 \pi} \W(\mu).
\end{equation}
In particular, if $\ \theta^2(\mu,\infty) < \infty\ $,
then $\ spt\ \mu\ $ has only finitely many components,
and hence these components are open and closed in $\ spt\ \mu\ $.

\end{enumerate}
\defin

%%%%%

\end{appendix}

%%%%%
%%%%%

%%%%%

\end{document}